\documentclass[final]{amsart}
\usepackage{pdfsync}
\usepackage{amsmath,amssymb}
\usepackage{enumerate}
\usepackage{xspace}
\usepackage{boxedminipage}
\usepackage{algorithmic}
\usepackage{listings}
\usepackage{tabularx}
\usepackage{subfigure}
\usepackage{tristan}
\usepackage{fullpage}

\renewcommand{\bi}[2]{\ensuremath{\cA\qp{#1,#2}}}
\renewcommand{\bih}[2]{\ensuremath{\cA_h\qp{#1,#2}}}

\renewcommand{\vec}[1]{\geovec{#1}}
\renewcommand{\mat}[1]{\geomat{#1}}

\renewcommand{\hoz}{\sobhz1(\W)}

\newcommand{\HCT}{\text{HCT}(k+2)}

\newcommand{\dx}{\,{\rm d}x}
\newcommand{\dS}{\,{\rm d}s}

\newcommand{\skel}{\Gamma}
\newcommand{\skelint}{\Gamma_{\rm int}}
\renewcommand{\avg}[1]{\{\kern-1.2mm\{#1\}\kern-1.2mm\}}

\renewcommand{\highlight}[1]{#1}

\numberwithin{equation}{section}
\setboolean{showtodo}{true}
\setboolean{shownotes}{true}
\setboolean{showchanges}{false}
\setboolean{usemathrsfs}{false}
\usepackage{tikz}
\usetikzlibrary{calc}

\author{
  Emmanuil H.~Georgoulis
}
\address{
  Emmanuil H.~Georgoulis,
  \thanks{
    Department of Mathematics,
    University of Leicester,
    University Road,
    Leicester LE1 7RH,
    UK and School of Applied Mathematical and Physical Sciences, National Technical University of Athens, Zografou 15780, Greece.
    {\tt{Emmanuil.Georgoulis@le.ac.uk}}
    }
}

\author{
  Charalambos G.~Makridakis
}
\address{
  Charalambos G.~Makridakis, 
  \thanks{Institute for Applied and Computational Mathematics-FORTH and DMAM University of Crete, Heraklion-Crete, Greece, GR 70013 and Department of Mathematics, University of Sussex, Brighton BN1 9QH, UK.
   \newline
    {\tt{C.Makridakis@sussex.ac.uk}}}
}

\author{
  Tristan Pryer
}
\address{
   Tristan Pryer,
   \thanks{
     Department of Mathematics and Statistics,
     University of Reading,
     Whiteknights,
     PO Box 220,
     Reading RG6 6AX,
     UK.
  {\tt{T.Pryer@reading.ac.uk}}
}}

\title[Babu\v ska-Osborn techniques in discontinuous Galerkin methods] 
{Babu\v ska-Osborn techniques in discontinuous Galerkin methods: \\ $\leb{2}$-norm error estimates for unstructured meshes}
\date{\today}
\pdfformat{true}

\begin{document}
\maketitle
\begin{abstract}
  We prove the inf-sup stability of the interior penalty class of
  discontinuous Galerkin schemes in unbalanced mesh-dependent norms,
  under a mesh condition allowing for a general class of meshes, which
  includes many examples of geometrically graded element
  neighbourhoods. \highlight{The inf-sup condition results in the
    stability of the interior penalty Ritz projection in a mesh
    dependent $\leb{2}$-norm, which allows for the proof of novel a priori error
    estimates that do \emph{not} depend on the global maximum meshsize in $\leb{2}$.
    Quasi-optimality results are also derived and some numerical
    experiments are given.}
\end{abstract}

\section{Introduction}

Discontinuous Galerkin (dG) methods are a popular family of
non-conforming finite element-type approximation schemes for partial
differential equations (PDEs) involving discontinuous approximation
spaces. In the context of elliptic problems their inception can be
traced back to the 1970s \cite{Nitsche:1971,Baker:1977,Arnold:1982}; see also
\cite{ArnoldBrezziCockburnMarini:2001} for an overview and
history of these methods for second order problems. 
For higher order problems, for example the (nonlinear) biharmonic problem, dG methods are a useful alternative to using $\cont{1}$-conforming elements  \cite{Baker:1977,SuliMozolevski:2007,GeorgoulisHouston:2009,GeorgoulisHoustonVirtanen:2011,Pryer:2014}, whose implementation (especially in the context of $hp$-version finite elements) can become complicated.

The derivation of $\leb{2}$-norm a priori error estimates is standard in the literature: for a standard dG method (e.g., symmetric interior penalty), for the Poisson problem with standard boundary conditions, and for piecewise linear finite elements, a combination of $\sobh{1}$ bounds and a duality approach yield the bound
\begin{equation}
  \Norm{u - u_h}_{\leb{2}(\W)}
  \leq
  C \max_{K\in\T{}} h_K
  \Big(\sum_{K\in\T{}} h_K^2 \Norm{D^2 u}_{\leb{2}(K)}^2\Big)^{1/2},
\end{equation}
i.e., the bound is identical to the respective bound for conforming finite element methods. It is well known that such a bound is not sharp: it is often desirable to use non-quasiuniform meshes generated, for instance, through an adaptive mesh refinement algorithm. In \cite{Makridakis:2016}, it was shown that this bound can be improved under some assumptions on the mesh to
\begin{equation}\label{result}
  \Norm{u - u_h}_{\leb{2}(\W)}
  \leq C
  \Big(\sum_{K\in\T{}} h_K^4 \Norm{D^2 u}_{\leb{2}(K)}^2\Big)^{1/2},
\end{equation}
for the conforming finite element method.
In this work, we prove \eqref{result} for the symmetric interior penalty dG method, thereby extending the results from \cite{Makridakis:2016} into the dG setting, under similar mesh assumptions. The mesh assumption, informally speaking, reads $\|\jump{ h}/\avg{h}\|_{\leb{\infty}(\skelint)}\le \alpha $, for some $0\le \alpha <1$, \highlight{sufficiently small} with $h$ denoting an element-wise constant function characterising the local meshsize and $\jump{\cdot}$ and $\avg{\cdot}$ the jump and average across the internal mesh skeleton $\skelint$. This effectively restricts the level of grading allowed on the underlying mesh, nonetheless allowing for geometrically graded meshes arising from adaptive mesh refinement procedures for example. 

The proof of \eqref{result} relies on a new inf-sup condition shown for unbalanced $\leb{2}$ and $\sobh{2}$-like mesh-dependent norms like those used in \cite{GeorgoulisPryer:2016}, however builds on this making use of new localisation techniques developed in \cite{Makridakis:2016} for the conforming finite element method and resolves a number of technical difficulties specific to the dG setting. In particular, in contrast to the conforming case, local bounds for the interface terms arising in the interior penalty dG bilinear form have to be also treated using non-standard ``bubble''-function techniques.
\highlight{At the same time, contrary to the respective conforming results in \cite{Makridakis:2016}, a new feature of our proof for the interior penalty dG method is that we do not use  the super-approximation arguments of \cite{NitscheSchatz:1974}. A side implication of our approach, is an improved dependence on the polynomial degree of the mesh restrictions, see below for details. }

This is in keeping with the spirit of the seminal work of Babu{\v s}ka and Osborn \cite{BabuskaOsborn:1980}, see also \cite{BabuskaOsbornPitkaranta:1980}, where the respective result to \eqref{result} for continuous finite element methods in one spatial dimension for second and fourth order problems are first proven. The present approach, however, is quite different on the technical level and results in inf-sup stability for $\leb{2}$- and $\sobh{2}$-like mesh-dependent norms under the aforementioned mesh assumption. Other potential applications of the analysis presented below include the development of convergent adaptive dG schemes for the $\leb{2}$-norm error, which would follow the respective developments of \cite{DemlowStevenson:2011} for conforming finite element methods,
quasi-best approximation results for nonconforming methods for elliptic \cite{VeeserZanotti:2017} and for evolution problems \cite{MakridakisBabuv-ska:1997}. 

We also take the opportunity to extend the ideas of \cite{Gudi:2010} into the $\leb{2}$ setting. This allows us to circumvent regularity restrictions that would require $u\in\sobh{s}$ for $s > 3/2$. Our analysis is quite general and holds for functions $u\in \sobh{1}$ only. The tools used to prove this result include an $\sobh{2}$-conforming reconstruction operator used in the a posteriori analysis of fourth order problems \cite{GeorgoulisHoustonVirtanen:2011}.

\section{Model problem and discretisation}
\label{sec:setup}

To assist the exposition of the key ideas, we shall consider the Poisson problem with homogeneous Dirichlet boundary conditions 
as model problem. The results presented in this work can be also proven with straightforward modifications for more general elliptic problems, such as ones with variable diffusivity and/or non-homogeneous boundary conditions.


More specifically, let $\W \subset \reals^d$ be an open Lipschitz domain and consider the problem: find $u\in \hoz$, such that
\begin{equation}
  \label{eq:weakform}
  \bi{u}{v} = \ltwop{f}{v} \Foreach v\in\hoz,
\end{equation}
where $\ltwop{\cdot}{\cdot}$ denotes the $\leb{2}$ inner product and
the bilinear form $\cA : \hoz \times \hoz \to \reals$ is given by
\begin{equation}
  \label{eq:bilinear-form}
  \bi{u}{v} := \int_\W \nabla u\cdot\nabla v \dx.
\end{equation}
Now, if $\Omega$ is such that $\Delta u\in \leb{2}(\W)$, we can also consider the \emph{unbalanced} bilinear form
$\cA : \sobh{2}(\T{})\cap \hoz \times \leb{2}(\W) \to \reals$ given by
\begin{equation}
  \label{eq:bilinear-form2}
  \bi{u}{v} := -\int_\W \Delta u v \dx,
\end{equation}
whose stability can be inferred via an inf-sup condition.
\begin{Pro}[inf-sup stability of the Laplacian]
  \label{pro:cont-inf-sup}
  With $\cA$ defined as in (\ref{eq:bilinear-form2}) we have that
  \begin{equation}
    \sup_{v\in\leb{2}(\W)} \frac{\bi{u}{v}}{\Norm{v}_{\leb{2}(\W)}}
    =
    \Norm{\Delta u}_{\leb{2}(\W)}.
  \end{equation}
  Also, assuming that $\Omega$ is convex, then the Miranda-Talenti inequality $\Norm{u}_{\sobh{2}(\W)}\le C \Norm{\Delta u}_{\leb{2}(\W)} $ holds for some $C>0$ independent of $u$, $f$ and we have the a priori bound
  \begin{equation}\label{pde_apriori}
    \norm{u}_{\sobh{2}(\W)} \leq C_{reg}\Norm{f}_{\leb{2}(\W)},
  \end{equation}
  for some $C_{reg}>0$ also independent of $u$ and of $f$.
\end{Pro}
\begin{Proof}
  The proof is immediate upon application of the Cauchy-Schwarz inequality on \eqref{eq:bilinear-form2}.
\end{Proof}

\subsection{Discretisation}
Let $\T{}$ be a conforming mesh of $\W\subset \reals^d$ into simplicial and/or box-type elements,
namely, $\T{}$ is a finite family of sets such that
\begin{enumerate}
\item $K\in\T{}$ implies $K$ is an open simplex (segment for $d=1$,
  triangle for $d=2$, tetrahedron for $d=3$) or an open box (quadrilateral for $d=2$, hexahedron for $d=3$),
\item for any $K,J\in\T{}$ we have that $\closure K\meet\closure J$ is either empty or
  a complete $(d-r)$-dimensional simplex/box (i.e., it is either a vertex for $r=d$, an
  edge for $r=d-1$, a face for $r=d-2$ when $d=3$, or the whole of $\closure K$ and $\closure J$) of both
  $\closure K$ and $\closure J$ and
\item $\union{K\in\T{}}\closure K=\closure\W$.
\end{enumerate}
The shape regularity constant of $\T{}$ is defined as
\begin{equation}
  \label{eqn:def:shape-regularity}
  \mu(\T{}) := \inf_{K\in\T{}} \frac{\rho_K}{h_K},
\end{equation}
where $\rho_K$ is the radius of the largest inscribed ball of
$K$ and $h_K$ is its diameter. An indexed family of
triangulations $\setof{\T n}_n$ is called \emph{shape regular} if 
\begin{equation}
  \label{eqn:def:family-shape-regularity}
  \mu:=\inf_n\mu(\T n)>0.
\end{equation}
%
For $s>0$, we define the \emph{broken} Sobolev space $\sobh{s}(\T{})$, by
\[
\sobh{s}(\T{}):= \{w\in \leb{2}(\W): w|_K\in \sobh{s}(K), K\in\T{}\},
\]
along with the broken gradient and Laplacian $\nabla_h\equiv \nabla_h(\T{})$  and $\Delta_h\equiv \Delta_h(\T{})$, i.e., the element-wise gradient and Laplacian operators.

We consider the \emph{finite element space}
\begin{gather}
  \label{eqn:def:finite-element-space}
  \fes :=   \{ \phi\in \leb{2}(\W):\phi|_K \in \poly k (K) \}
\end{gather} where $\poly k (K)$ is the space of polynomials of total degree $k$ \highlight{for $k \geq 1$}. Alternatively, when $K\in\T{}$ is a box-type element, we can also consider polynomials of degree $k$ in each variable, typically mapped from a reference hypercube. {\highlight{Apart from assuming shape-regularity for the remainder of this work}, the fine properties of the respective finite element spaces are of no essential consequence to the results below, as long as standard best approximation bounds are available for the elements considered. 

Let also $\skel=\cup_{K\in\T{}}\partial K$ denote the skeleton of the
mesh $\T{}$ and set $\skelint:=\skel\backslash\partial\W$ to denote the skeleton interior to $\W$. 

\begin{Defn}[jumps and averages]
  \label{defn:averages-and-jumps}
  We define average and jump operators for arbitrary scalar
  $v\in \sobh{s}(\T{})$ and vector $\vec v\in [\sobh{s}(\T{})]^d$ functions, with $s>3/2$,  as
  \begin{gather}
    \label{eqn:average}
    \avg{v} =  {\frac{1}{2}\qp{v|_{K_1} + v|_{K_2}}},
    \qquad \avg{\vec v} = {\frac{1}{2}\qp{\vec{v}|_{K_1} + \vec{v}|_{K_2}}},
    \\\nonumber\\
    \label{eqn:jump}
    \jump{v} = {{{v}|_{K_1} \geovec n_{K_1} + {v}|_{K_2}} \geovec n_{K_2}},
    \qquad
    \jump{\vec v}
    = 
    \Transpose{\qp{\vec{v}|_{K_1}}}\geovec n_{K_1} + \Transpose{\qp{\vec{v}|_{K_2}}}\geovec n_{K_2}.
  \end{gather}
  Note that on the boundary of the domain $\partial\W$ the jump and
  average operators are defined as
  \begin{gather}
    \avg{v}
    \Big\vert_{\partial\W} 
    := v,
    \qquad 
    \avg{\geovec v}
    \Big\vert_{\partial\W}
    :=
    \geovec v,
    \\
    \jump{v}
    \Big\vert_{\partial\W}
    := 
    v\geovec n,
    \qquad 
    \jump{\geovec v}
    \Big\vert_{\partial\W} 
    :=
    \Transpose{\geovec v}\geovec n,
  \end{gather}
\end{Defn}

Further, we define $\funk h\W\reals_+$ to be the {piecewise
constant} \emph{meshsize function} of $\T{}$ given by $h |_K:=h_K$, $K\in\T{}$ and $h|_\skel:=\avg{h}$. The conformity assumption of the mesh, along with shape regularity imply the equivalence
\begin{equation}
  \label{eq:qu-const}
  C_{qu}^{-1}h_K\le h(x) \le C_{qu}h_K, 
\end{equation}  
for all $x\in \omega_K:=\cup_{K'\in\T{}:\bar{K}\cap \bar{K}'\neq \emptyset}K'$, for some $C_{qu}>0$ depending only on $\mu$.

\subsection{Interior penalty discontinuous Galerkin method}
We consider the interior penalty (IP) discontinuous Galerkin discretisation of
(\ref{eq:bilinear-form}), reading: find $u_h \in \fes$ such that
\begin{equation}\label{eq:dg}
  \bih{u_h}{v_h} = \ltwop{f}{v_h} \Foreach v_h \in \fes,
\end{equation}
where  
\begin{equation}
  \label{eq:IP}
  \begin{split}
    \bih{u_h}{v_h} 
    &=
    \int_\W \nabla_h u_h \cdot \nabla_h v_h \dx
    -
    \int_{\skel} \big(\jump{v_h} \cdot \avg{\nabla u_h} +\theta\jump{u_h} \cdot \avg{\nabla v_h} 
    - \sigma \jump{u_h}\cdot \jump{v_h}\big)\dS
    ,
  \end{split}
\end{equation}
where $\sigma > 0$ is the, so-called, \emph{discontinuity penalisation parameter} given by 
\begin{equation}
  \label{eq:sigma}
  \sigma:= C_{\sigma}\frac{k^2}{h},
\end{equation}
and $\theta\in [-1,1]$ a (global) constant used to select between the symmetric IP dG method ($\theta=1$) and its non-symmetric variant ($\theta=-1$). As expected optimal results are obtained when $\theta = 1$. For completeness we also discuss the case of $\theta \neq 1$ (see Remark \ref{rem:ns-int-pen}). The constant $C_{\sigma}>0$ is also typically chosen globally: when $\theta=1$ it should be chosen large enough so as to counteract a constant of an inverse estimate to achieve coercivity, while it can be chosen freely when $\theta=-1$. Numerical evidence suggests that the choice $\theta=-1$ results in dG methods which converge suboptimally with respect to the meshsize $h$ for even polynomial degrees, when the error is measured in the $\leb{2}$-norm \cite{HSS02,Harriman_et_al,Hartmann}.

\begin{Defn}[mesh dependent norms]
  \label{def:mesh-dep-norms}
  We introduce the {mesh dependent} $\leb{2}$, $\sobh1$ and $\sobh2$ norms to be 
  \begin{gather}
    \znorm{w}^2 := \Norm{w}_{\leb{2}(\W)}^2
    +
    \Norm{h^{3/2} \avg{\nabla w}}_{\leb{2}(\skel)}^2
    +
    \Norm{h^{1/2} \avg{w}}_{\leb{2}(\skelint)}^2
    +\Norm{h^{1/2} \jump{w}}_{\leb{2}(\skel)}^2
    \\
    \enorm{w}^2 := \Norm{\nabla_h w}_{\leb{2}(\W)}^2 + \Norm{\sqrt{\sigma}\jump{w}}_{\leb{2}(\skel)}^2
    \\
    \eenorm{w}^2 
    :=
   \Norm{\nabla_h w}_{\leb{2}(\W)}^2 
    +
    \Norm{\Delta_h w}_{\leb{2}(\W)}^2
    +
    \Norm{h^{-1/2}\jump{\nabla w}}_{\leb{2}(\skelint)}^2
    +
    \Norm{ h^{-3/2}\jump{w}}_{\leb{2}(\skel)}^2.
  \end{gather} 
\end{Defn}

\begin{Rem}[motivation and properties of mesh dependent norms]
  The motivation for the norms given in Definition \ref{def:mesh-dep-norms} is that, upon integration by parts, the IP dG bilinear form becomes
  \begin{equation}
    \label{eq:IP2}
    \begin{split}
      \bih{u_h}{v_h} 
      &=
      -\int_\W \Delta_h u_h \ v_h \dx
      +
      \int_{\skelint} \jump{\nabla u_h} \avg{v_h}\dS 
      - \int_{\skel} \theta\jump{u_h} \cdot \avg{\nabla v_h}
      -
      \sigma \jump{u_h}\cdot \jump{v_h}\dS,
    \end{split}
  \end{equation}
  whence, for $w,v\in \sobh{2}(\T{})$, 
  \begin{equation}
    \label{eq:cont-bound}
    \norm{\bih{w}{v}} \leq C\eenorm{w} \znorm{v}.
  \end{equation}
  Notice that the norm $\eenorm{\cdot}$ includes $\enorm{\cdot}$ to ensure it is, indeed, a norm.
  
  The norm $\znorm{\cdot}$ is also equivalent to the $\leb{2}$ norm over $\fes$ in view of standard inverse inequalities, that is, for any $w_h\in\fes$ there exists a $C>0$ such that
  \begin{equation}
    \label{eq:norm-equiv}
    C^{-1}\znorm{w_h} \leq \Norm{w_h}_{\leb{2}(\W)} \leq \znorm{w_h}.
  \end{equation}
\end{Rem}

\begin{Pro}[continuity and coercivity of $\bih{\cdot}{\cdot}$ in $\enorm{\cdot}$]
  \label{pro:cont-coer}
  For $C_\sigma$ large enough, the bilinear form $\bih{\cdot}{\cdot}$ satisfies
  \begin{gather}
    \bih{u_h}{u_h} \geq c_0 \enorm{u_h}^2\label{eq:coer}
    \\
    \bih{u_h}{v_h} \leq C_0 \enorm{u_h} \enorm{v_h},
  \end{gather}
  for $c_0,C_0>0$ independent of $h$, $C_{\sigma}$, $u_h$, and $v_h$.
  \end{Pro}
Lax-Milgram Theorem guarantees a unique solution to the problem (\ref{eq:IP}).
%

The main result of this work is the following theorem, the proof of which we shall dedicate Section~\ref{sec:proof} to.
\begin{The}[Inf-sup stability of the dG method]
  \label{the:inf-sup}
  Let $\bih{\cdot}{\cdot}$ be the bilinear form given in (\ref{eq:IP}) with $\theta=1$, and assume that the penalty parameter $\sigma$ is chosen large enough to ensure the validity of \eqref{eq:coer}. Suppose that the underlying mesh of the finite element space satisfies 
  \begin{equation}\label{mesh_ass}
  \|\jump{ h}/\avg{h}\|_{\leb{\infty}(\skelint)}\le \alpha \ ,\quad  \text{ for some } 0\le \alpha <1\ \text{small enough,}
  \end{equation}  
  \highlight{depending on the shape-regularity constant of the mesh and on the polynomial degree $k$.}
   Then, there exists a constant $\gamma>0$, independent of $h$, $w_h$ and $v_h$, such that
  \begin{equation}
    \sup_{0\neq v_h\in\fes} \frac{\bih{w_h}{v_h}}{\eenorm{v_h}} \geq \gamma \znorm{w_h} \Foreach w_h\in\fes.
x  \end{equation}
\end{The}

\begin{Rem}[Construction of meshes satisfying $\|\jump{ h}/\avg{h}\|_{\leb{\infty}(\skelint)}\le \alpha$ for any $0<\alpha<1$]
  \label{Rem:construction}
  To show that the mesh condition, $\|\jump{
    h}/\avg{h}\|_{\leb{\infty}(\skel)}\le \alpha$, is not restrictive
  and still allows for highly graded meshes we illustrate the
  construction of a highly non-uniform mesh satisfying this condition. In particular, as we shall see below, \emph{geometrically graded} meshes are admissible.

Let $0< \beta<1$ and $\Omega=(0,1)^2$. Next, consider a grid on
$\Omega=(0,1)^2$ given by the points
$\{0,\beta^N,\beta^{N-1},\dots,\beta,1\}$ in each direction, and
construct the respective structured rectangular mesh with elements
\[
\widetilde{K}_{ij}:=(\beta^i,\beta^{i+1})\times (\beta^j,\beta^{j+1}),\quad i,j=\{0,1,\dots, N\}\cup\{\infty\},
\]
making use of the convention $\beta^\infty=0$.  Let
$\mathcal{T}_{\beta}$ denote the triangular mesh constructed from the
subdivision $\{\tilde{K}_{ij}\}$ by taking the southwest-northeast
diagonal on each $\tilde{K}_{ij}$. We refer to Figure \ref{fig:construction} for an
illustration.
 
Let now $h_{ij}:=(1-\beta)\sqrt{\beta^{2i}+\beta^{2j}}$ be the
diameter of each of the two elements arising by taking the diagonal of
$\widetilde{K}_{ij}$. We begin by noting that $\jump{h}=0$ on all
``diagonal'', i.e., non-axiparallel internal faces of the mesh. Next,
upon observing that, adjacent elements to the two elements in
$\tilde{K}_{ij}$ have diameters $h_{(i+1)j},
h_{(i-1)j},h_{i(j+1)},h_{i(j+1)}$, respectively, we can continue the
argument for one case, say $h_{(i+1)j}$ only, without loss of
generality due to symmetry. In this case, on the common face between
these two elements, noting that $\jump{h^2} = 2\jump{h}\avg{h}$, we
have, respectively,
\begin{equation}\label{bound_grading}
\frac{\jump{ h}}{\avg{h}}=\frac{\jump{ h^2}}{2\avg{h}^2}\le \frac{2\jump{ h^2}}{\avg{h^2}}= \frac{4(1-\beta)^2( \beta^{2i}-\beta^{2(i+1)})}{(1-\beta)^2(\beta^{2i}+\beta^{2(i+1)}+2\beta^{2j})}= \frac{4(1-\beta^2)}{1+\beta^2+2\beta^{2(j-i)}},
\end{equation}
giving $\|\jump{ h}/\avg{h}\|_{\leb{\infty}(\skelint)}\to 0$ as $\beta\to 1^-$.

\end{Rem}

\newcommand*\rows{20}
\newcommand*\be{0.9}
\begin{figure}
  \begin{tikzpicture}[scale=6]
    \draw ($(1, \be^\rows*\be)$) -- ($(\be^\rows*\be,\be^\rows*\be)$);
    \draw ($(\be^\rows*\be, \be^\rows*\be)$) -- ($(\be^\rows*\be,1)$);
    \foreach \rowi in {0, 1, ...,\rows}
             {
               \foreach \rowj in {0, 1, ...,\rows}
                        {
                          \draw ($(\be^\rowi, \be^\rowj)$) -- ($(\be^\rowi*\be,\be^\rowj)$);
                          \draw ($(\be^\rowi, \be^\rowj)$) -- ($(\be^\rowi,\be^\rowj*\be)$);
                          \draw ($(\be^\rowi*\be, \be^\rowj*\be)$) -- ($(\be^\rowi,\be^\rowj)$);
                        }
             }
             \draw ($(\be^1, \be^1.5)$) -- ($(\be^2,\be^1.5)$)  node[below,right]{$K$};
             \filldraw[fill=blue!40!white, draw=black] ($(\be^1, \be^1)$) rectangle ($(\be^2,\be^2)$);
             \node at ($(\be^1.5, \be^1.5)$)   (a) {$\tilde{K}_{ij}$};
             \pgfmathparse{\rows};
  \end{tikzpicture}
  \caption{\label{fig:construction} Illustration of the graded mesh constructed in Remark
    \ref{Rem:construction}. Here $\beta = \be$, $N = \rows$, with \eqref{bound_grading} giving
    $\|\jump{ h}/\avg{h}\|_{\leb{\infty}(\skelint)} \le 0.42$.  }
\end{figure}

\begin{Rem}[Interpreting the condition $\|\jump{ h}/\avg{h}\|_{\leb{\infty}(\skelint)}\le \alpha$]
  Figure \ref{fig:mesh} shows three different classes of mesh, one being generated through a newest vertex bisection adaptive refinement procedure another being an artificially graded mesh and a third being of Shishkin type. In all cases the values of $\|\jump{ h}/\avg{h}\|_{\leb{\infty}(\skel)}$ are computed. As expected, standard, shape-regular locally adapted meshes generated through newest vertex bisection refinement satisfy $\Norm{\jump{h}/\avg{h}}_{\leb{\infty}(\Gamma)} \lesssim 1$. 
  
  \highlight{ Note that the mesh function we make use of in this work
    is different than that used in the works of
    \cite{Eriksson:1994,DemlowStevenson:2011,Makridakis:2016} they
    are, however, related. Let $\widetilde h$ denote the piecewise
    linear continuous mesh function defined in
    \cite{Eriksson:1994,DemlowStevenson:2011,Makridakis:2016}. It can
    be shown that the mesh function $\widetilde h$ coincides with the
    nodal averaged reconstruction of the piecewise constant
    discontinuous mesh function $h$ \cite{KP}; see also
    \cite{GeorgoulisPryer:2018} for some related ideas. It is, thus,
    possible to construct stability bounds to relate the two mesh
    functions. Indeed, modifying the argument of \cite[Lemma
      4.2]{DemlowGeorgoulis:2012} we have, 
    \begin{equation}\label{connection}
      \|\nabla \widetilde h\|_{\leb{\infty}(\W)}
      \leq
      C \Norm{ \jump{h}/\avg{h}}_{\leb{\infty}(\skelint)},
    \end{equation}
    where $C$ depends upon the shape regularity of the mesh only.  We
    also refer to Remark \ref{polydeg} below for a different aspect in
    the comparison between the two conditions.}

\end{Rem}

\begin{Rem}
  The main result of this work, Theorem \ref{the:inf-sup} is also
  valid if we replace \eqref{mesh_ass} by classical mesh condition
  from \cite{Eriksson:1994,DemlowStevenson:2011,Makridakis:2016} via
  the use of superapproximation results. We refrain from doing so in
  this work as \eqref{mesh_ass} appears to be more natural in this
  context of discontinuous finite element spaces, (cf. also the
  discussion in Remark \ref{mess_ass_2} below,) and we refer to
  \cite{Makridakis:2016} for the proof of the respective result for
  the conforming finite element method.
\end{Rem}

\begin{figure}[!h]
 \caption[]
         {\label{fig:mesh}
           Interpreting the condition $\|\jump{ h}/\avg{h}\|_{\leb{\infty}(\skelint)}\le \alpha$
 }
 \begin{center}
   \subfigure[{\label{fig:a1}
       Shape-regular adapted mesh;
       here $\Norm{\jump{h}/\avg{h}}_{\leb{\infty}(\Gamma)} \approx 0.2$.
   }]{
     \includegraphics[scale=\figscale,width=0.28\figwidth]{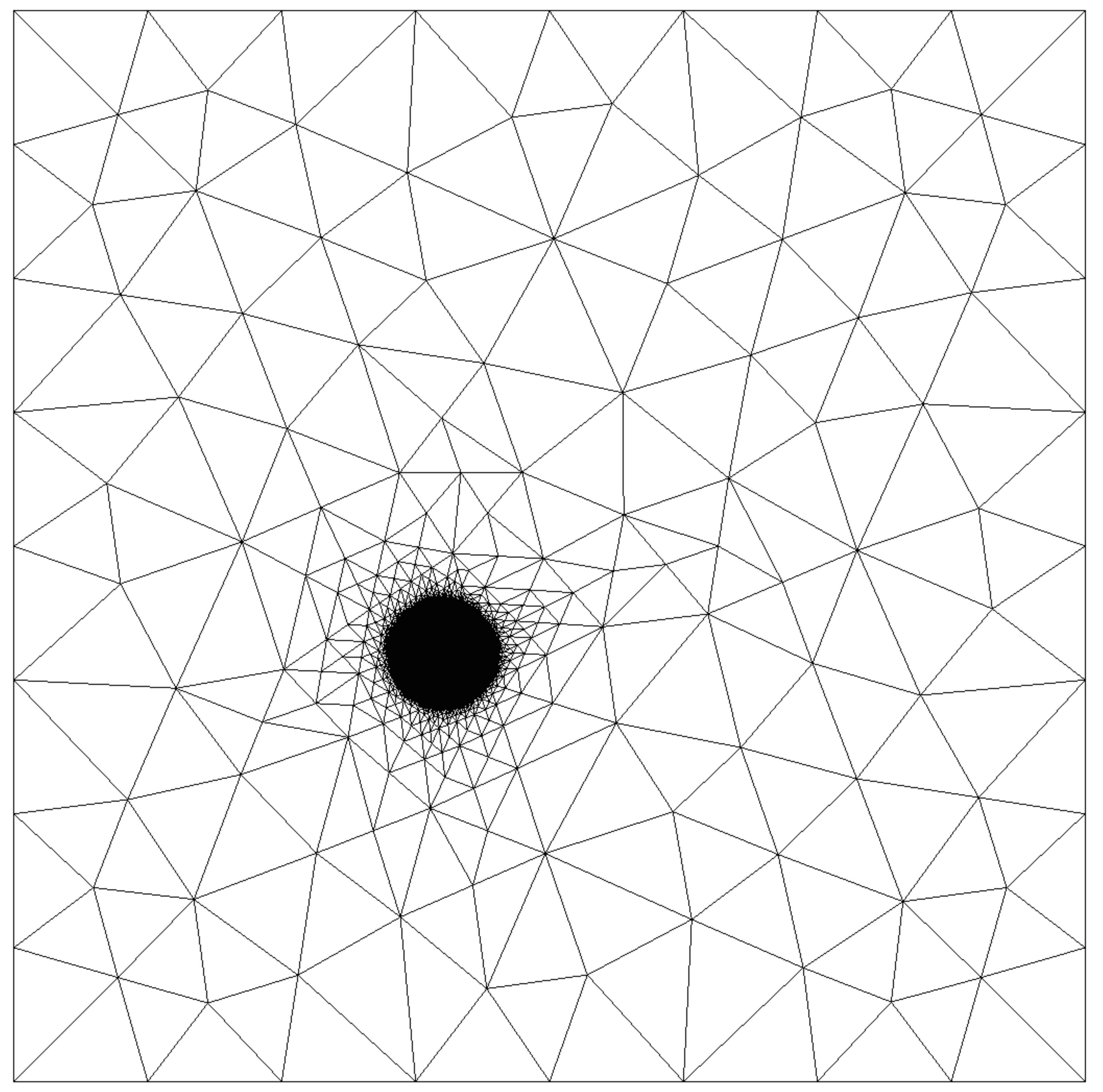}
   }\hspace{.2cm}
   \subfigure[{\label{fig:a2}
       An anisotropic graded mesh; here $\Norm{\jump{h}/\avg{h}}_{\leb{\infty}(\Gamma)} \approx 0.61$.
   }]{
     \includegraphics[scale=\figscale,width=0.28\figwidth]{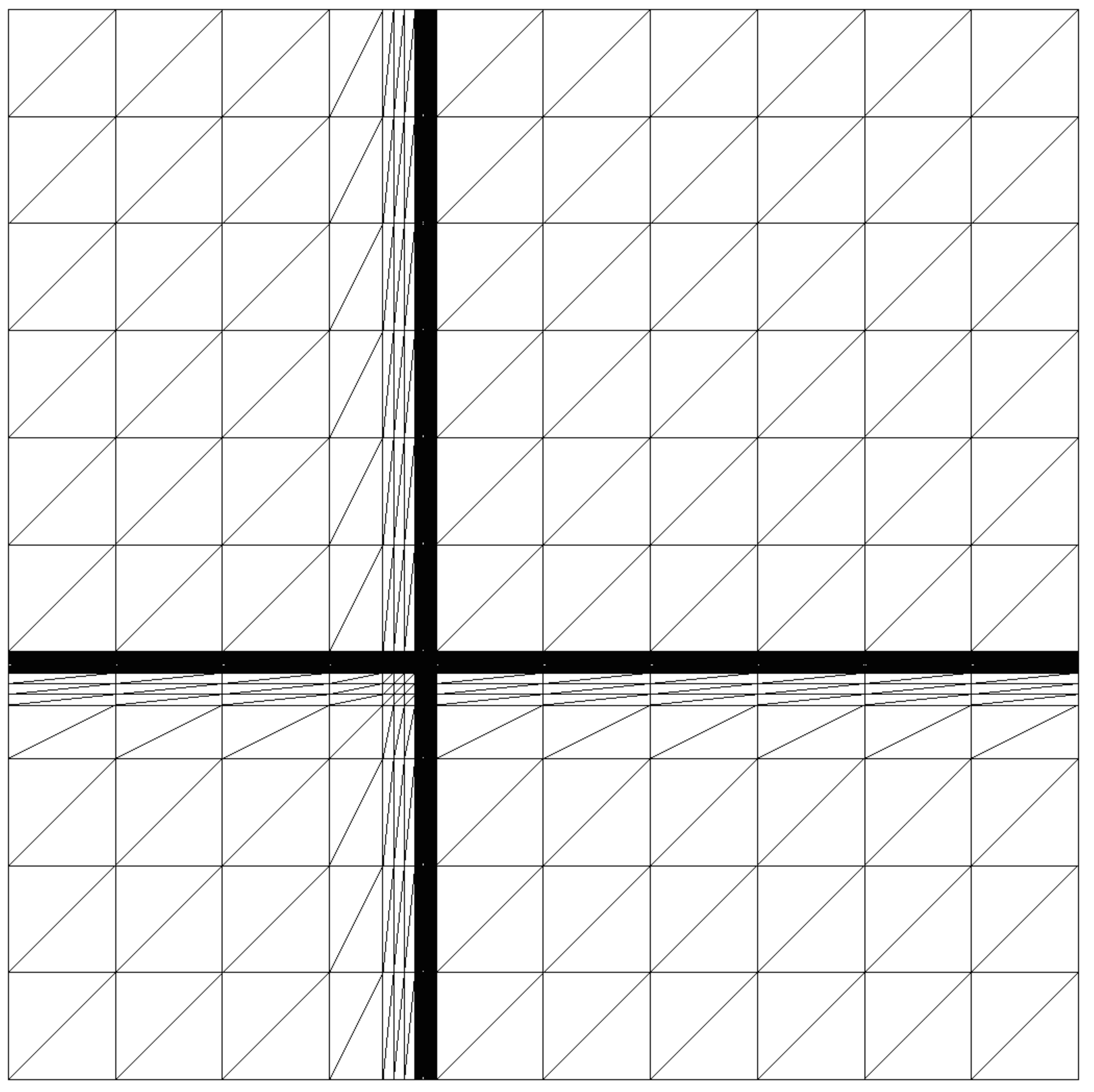}
   }\hspace{.2cm}
   \subfigure[{\label{fig:a2}
       A Shishkin-type mesh; here $\Norm{\jump{h}/\avg{h}}_{\leb{\infty}(\Gamma)} \approx 0.59$.
   }]{
     \includegraphics[scale=\figscale,width=0.28\figwidth]{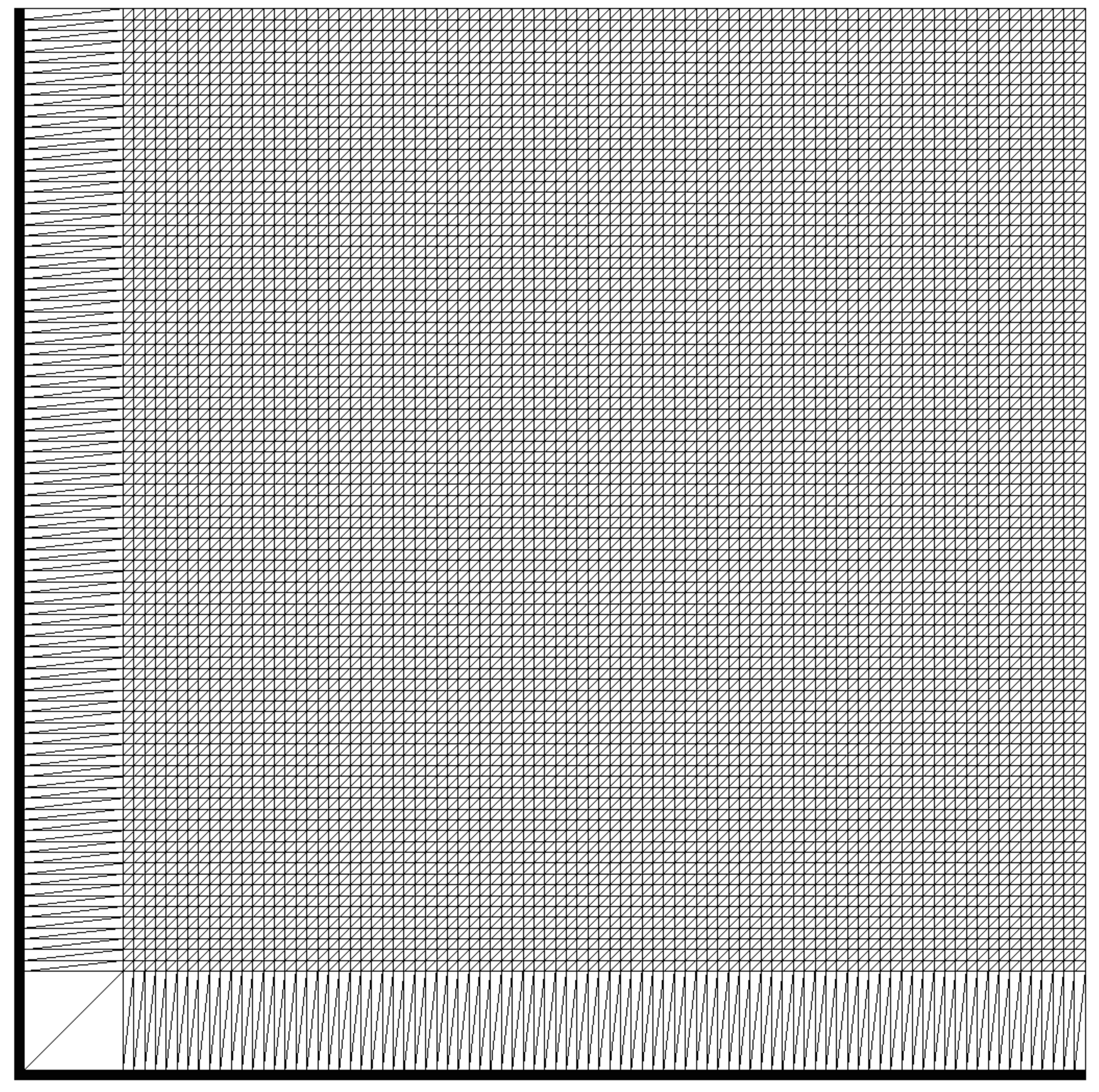}
   }
   \end{center}
 \end{figure}

Equipped with Theorem \ref{the:inf-sup}, we have the following result, stating the $\leb{2}$-norm error optimality of the interior penalty dG method under the mesh assumption \eqref{mesh_ass}. 

\begin{Cor}
  Let $R : \sobh2(\T{})\to\fes$ be the dG-Ritz-projection operator defined for $u\in\sobh{2}(\T{})$ by
  \begin{equation}
    \label{eq:ip-proj}
    \bih{R u}{v_h} = \bih{u}{v_h} \Foreach v_h \in \fes.
  \end{equation}
  Then under the hypotheses of Theorem \ref{the:inf-sup}, we have:
  \begin{enumerate}
  \item
    \label{it:1}
    $R$ is stable in $\znorm{\cdot}$, that is:
   $
      \znorm{Ru} \leq \gamma^{-1} \znorm{u}.
    $
     \item
    \label{it:2}
    $R$ satisfies quasi-optimal error bounds in $\znorm{\cdot}$, that is:
    \begin{equation}
      \znorm{u - Ru}
      \leq
      \qp{1+\gamma^{-1}} \inf_{v_h \in \fes} \znorm{u - v_h}.
    \end{equation}
  \item
    \label{it:3}
    If $u\in\sobh{k+1}(\W)$ solves \eqref{eq:weakform} and $u_h\in \fes$ solves (\ref{eq:dg}), then
    \begin{equation}
      \znorm{u - u_h}
      \leq
      C \sum_{K\in\T{}} \qp{\Norm{h^{k+1} D^{k+1} u}_{\leb{2}(K)}^2}^{1/2}.
    \end{equation}
  \end{enumerate}
\end{Cor}
\begin{Proof}
  For \eqref{it:1}, Theorem \ref{the:inf-sup}, the definition of $R$ (\ref{eq:ip-proj}) and the continuity bound (\ref{eq:cont-bound}), imply   \begin{equation}
    \begin{split}
      \gamma \znorm{R u}
      &\leq
      \sup_{0\neq v_h\in\fes}
      \frac{\bih{R u}{v_h}}{\eenorm{v_h}}
      \leq
      \sup_{0\neq v_h\in\fes}
      \frac{\bih{u}{v_h}}{\eenorm{v_h}}
    \leq
      \znorm{u}.
    \end{split}
  \end{equation}
  For \eqref{it:2}, note that for any $v_h \in \fes$
  \begin{equation}
    \begin{split}
      \znorm{u - R u}
      &\leq
      \znorm{u - v_h}
      +
      \znorm{R \qp{v_h - u}}
      \leq
      \qp{1+\gamma^{-1}}
      \znorm{u - v_h},
    \end{split}
  \end{equation}
  due to \eqref{it:1}. Finally, \eqref{it:3} follows by choosing $v_h$ to be an appropriate interpolant and using respective best approximation bounds.
\end{Proof}

\begin{figure}[!ht]
 \caption[]
         {\label{fig:sol}
           In this experiment we test the $\leb{2}$ convergence of the
           interior penality method and demonstrate that even for the
           worst class of mesh given in Figure \ref{fig:mesh}, optimal
           $\leb{2}$ convergence is achieved. Here we chose $C_\sigma
           = 20$, smaller values of $C_\sigma$ resulted in a
           suboptimal convergence in $\leb{2}$ norm.
         }
         \begin{center}
           \subfigure[{\label{fig:b1}
               After $1$ global refinement.
   }]{
     \includegraphics[scale=\figscale,width=0.47\figwidth]{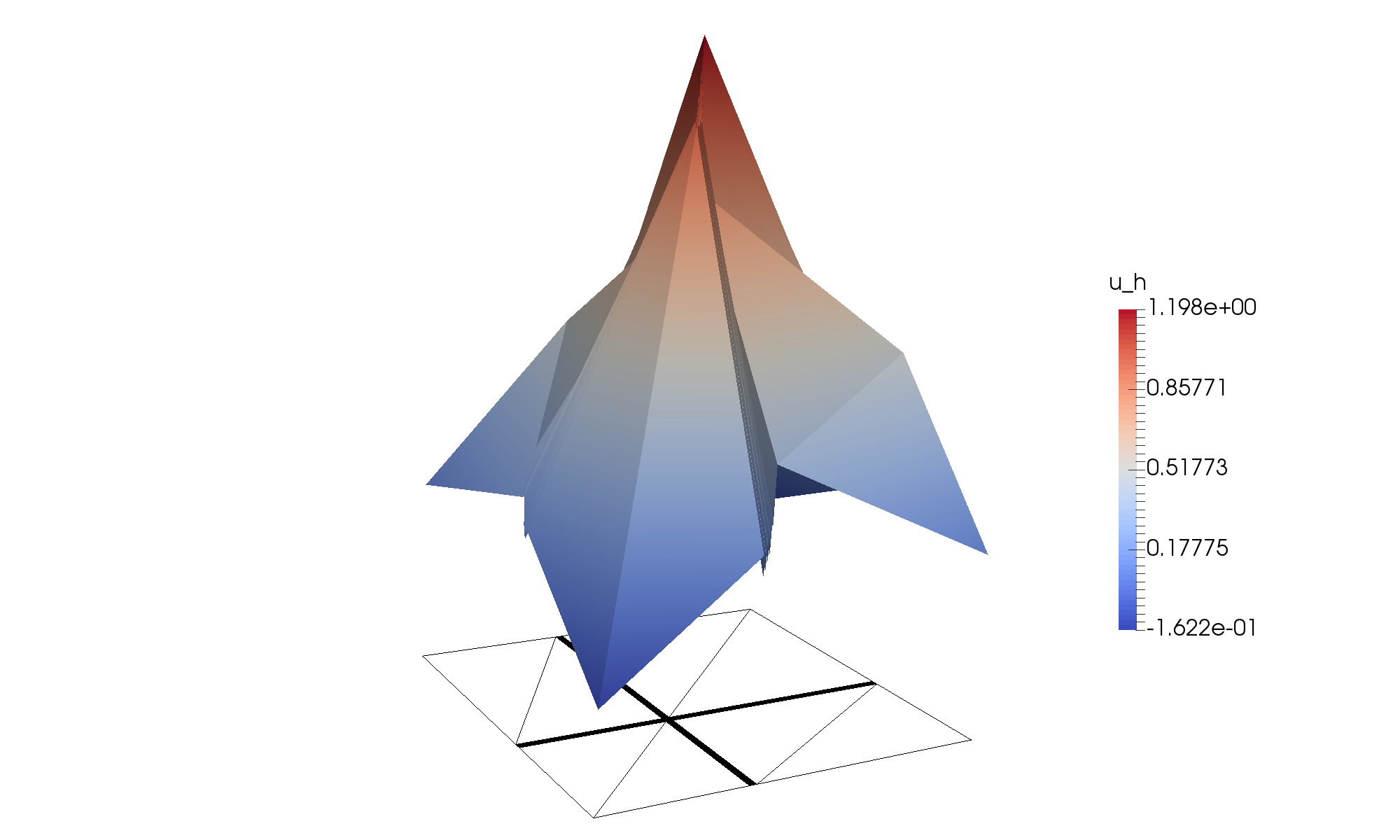}
   }
   \hfill
   \subfigure[{\label{fig:b2}
       After $2$ global refinements.
   }]{
     \includegraphics[scale=\figscale,width=0.47\figwidth]{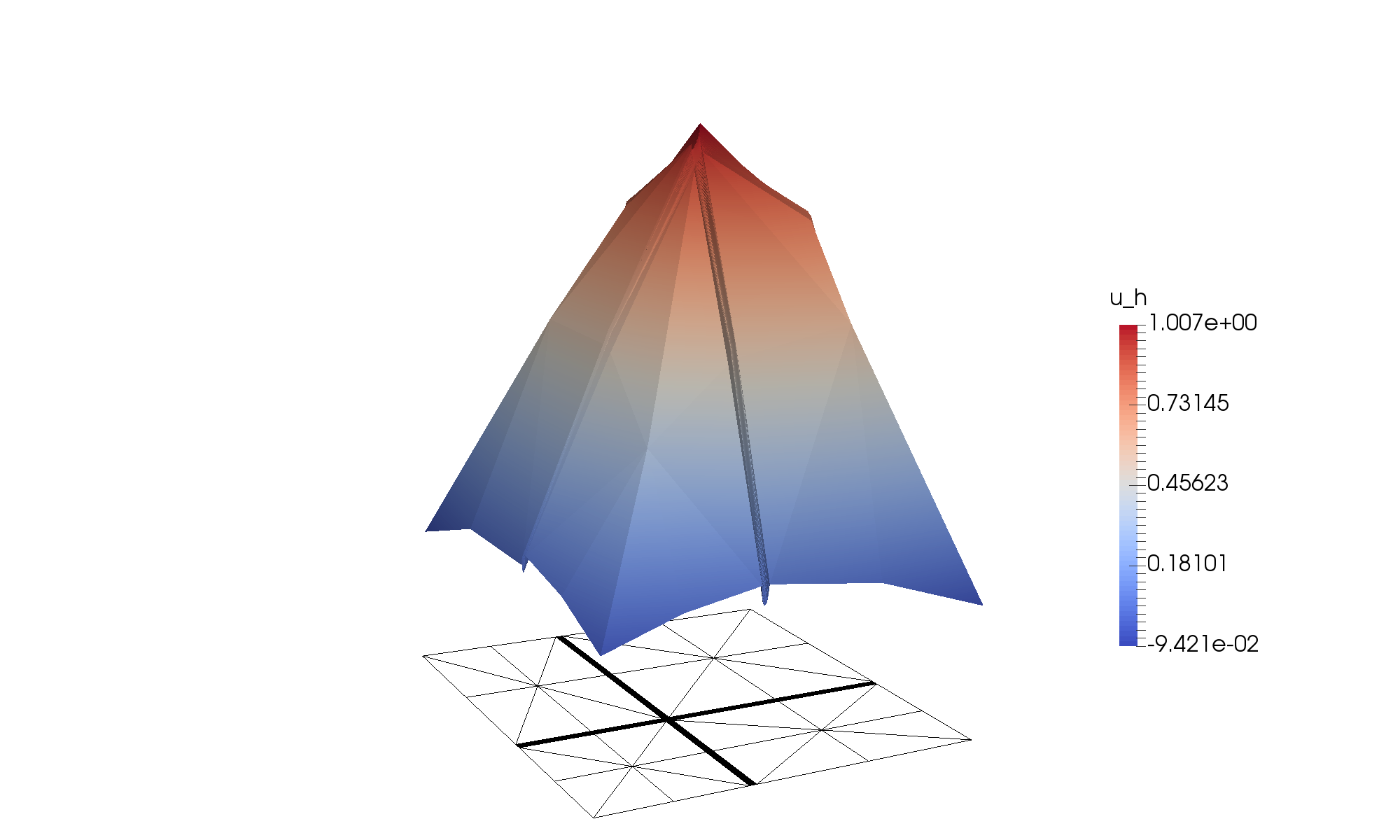}
   }
   \subfigure[{\label{fig:b3}
       After $3$ global refinements.
   }]{
     \includegraphics[scale=\figscale,width=0.47\figwidth]{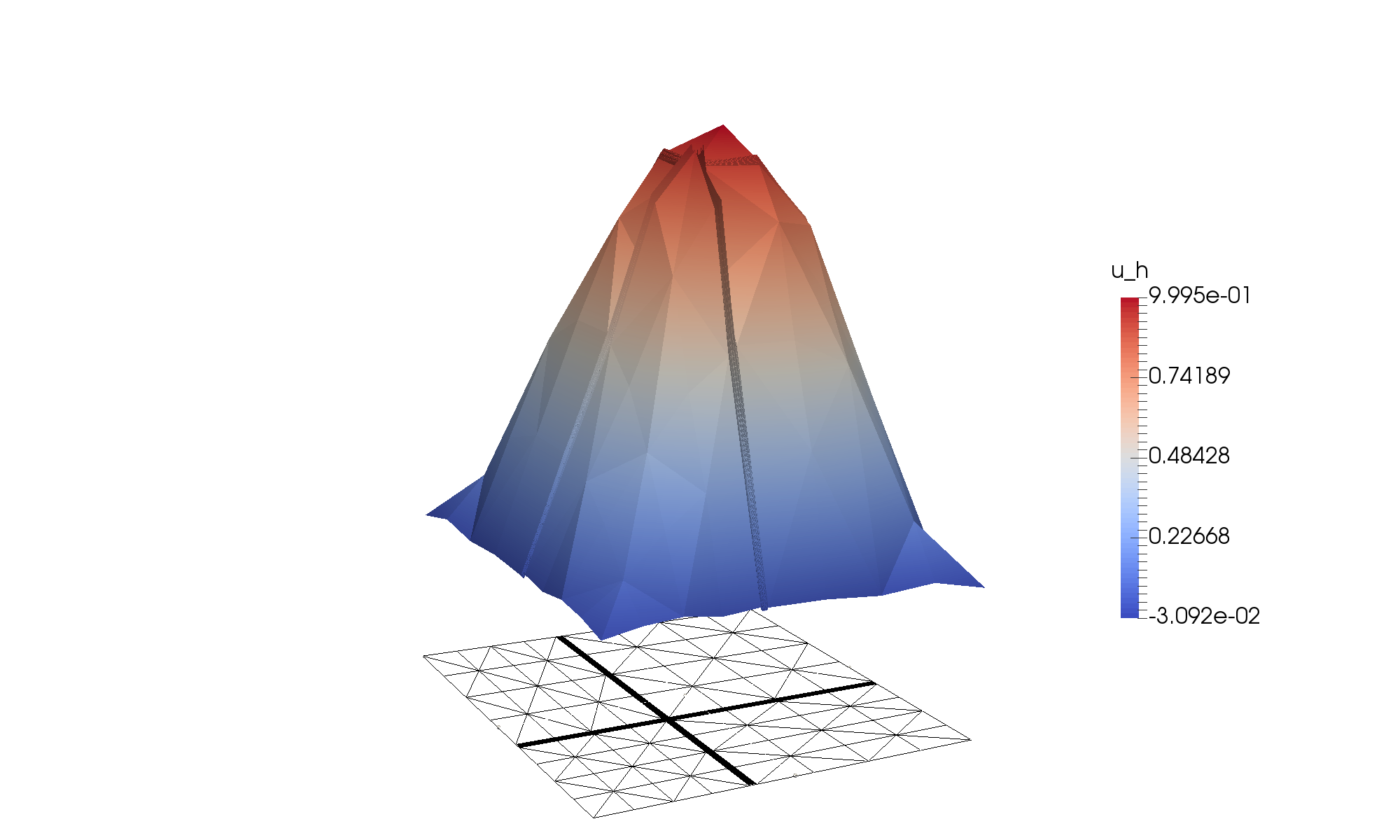}
   }
   \hfill
   \subfigure[{\label{fig:b4}
       After $4$ global refinements.
   }]{
     \includegraphics[scale=\figscale,width=0.47\figwidth]{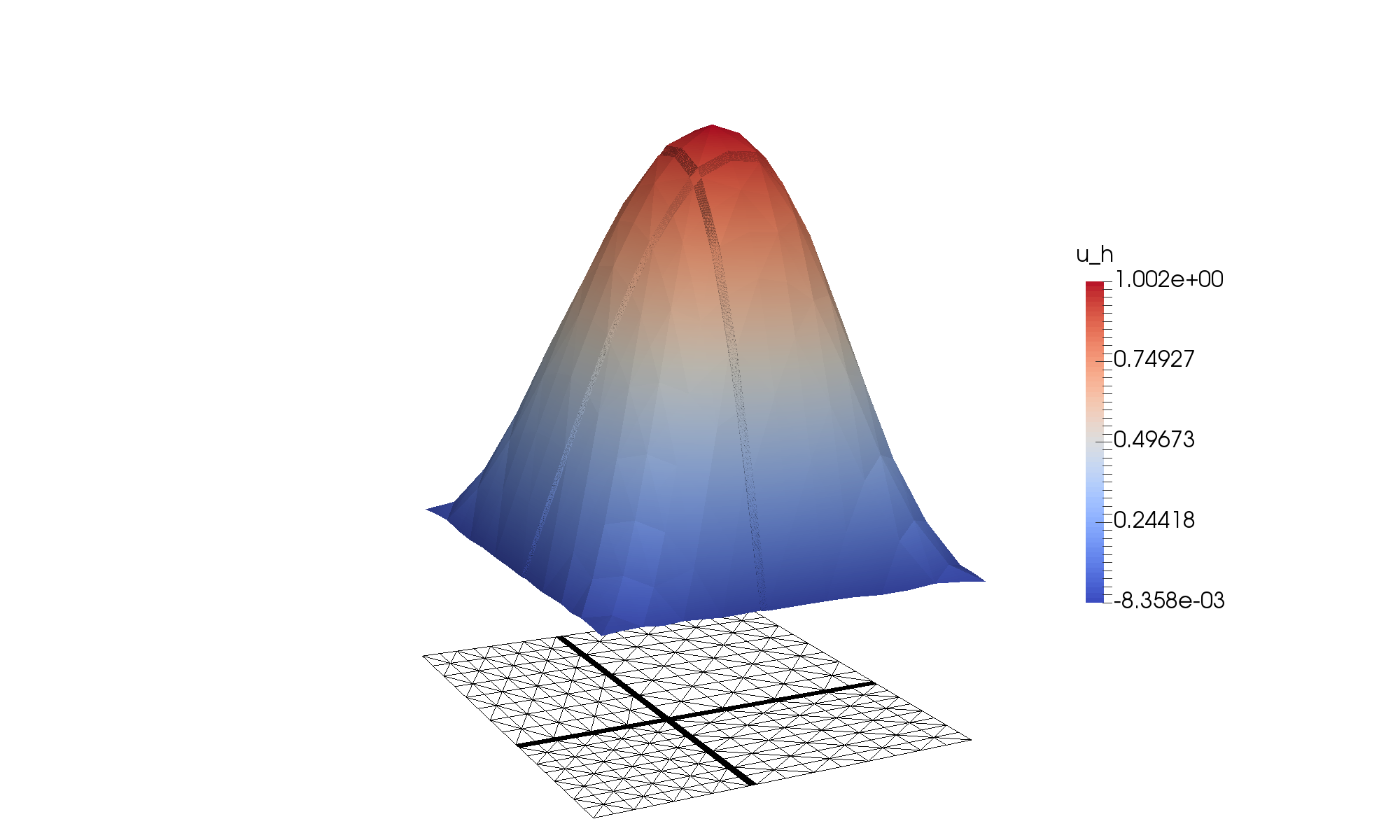}
   }
   \subfigure[{\label{fig:b3}
       After $5$ global refinements.
   }]{
     \includegraphics[scale=\figscale,width=0.47\figwidth]{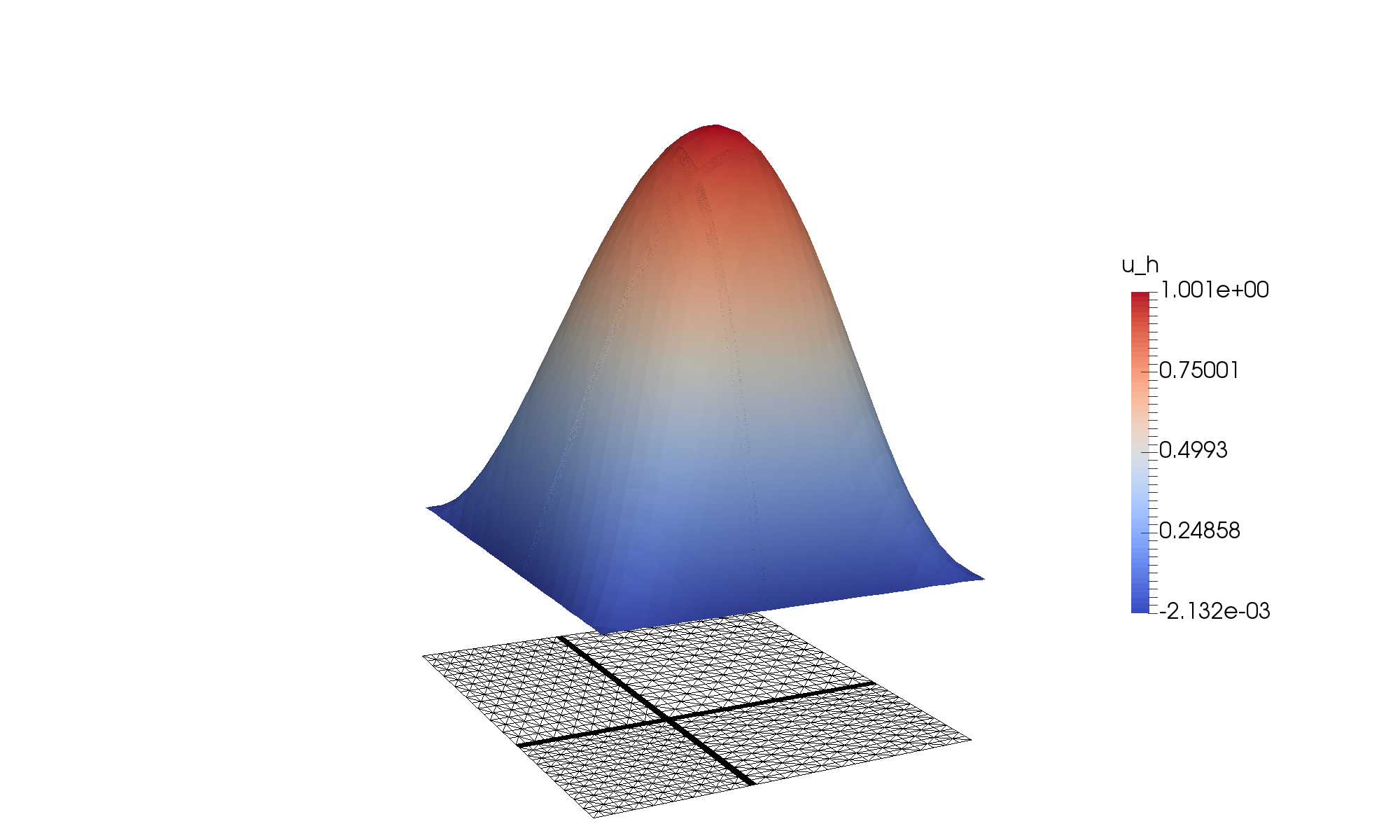}
   }
   \hfill
   \subfigure[{\label{fig:b4}
       Convergence rates of the approximation.
   }]{
     \includegraphics[scale=\figscale,width=0.47\figwidth]{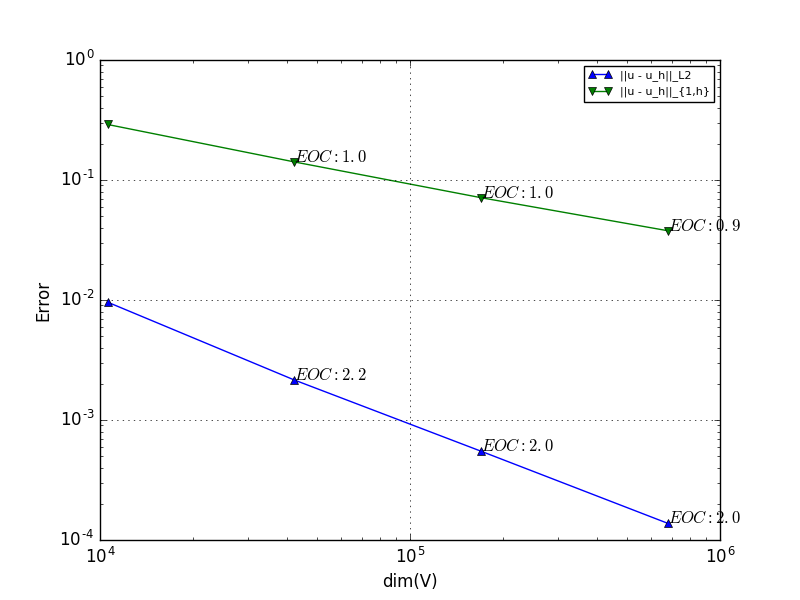}
   }
   \end{center}
 \end{figure}

\section{Proof of Theorem \ref{the:inf-sup}}
\label{sec:proof}


We begin by proving a crucial technical result regarding the stability of the dG-Ritz-projection operator in the $\leb{2}$-norm.

\begin{Lem}[$\leb{2}$-stability of $R$]
  \label{lem:bound-for-R}
Let $w\in\sobh{2}(\W)$ and assume that the mesh satisfies \eqref{mesh_ass}. Then, for $\theta=1$, its dG-Ritz-projection $Rw$ satisfies the bound
  \begin{equation}\label{key}
    \Norm{R w}_{\leb{2}(\W)}
    \leq 
    C \big(
      \Norm{h \nabla w}_{\leb{2}(\W)} + \Norm{w}_{\leb{2}(\W)} +  \Norm{h^{3/2}\avg{\nabla w}}_{\leb{2}(\skel)}  
    \big).
  \end{equation}
\end{Lem}
\begin{Proof}
Let $g\in\sobh{2}(\W)$ be the solution to 
  \begin{equation}
    \begin{split}
      -\Delta g &= R w \text{ in } \W,
      \quad
      g = 0 \text{ on } \partial\W,
    \end{split}
  \end{equation}
  for which we assume the a priori bound \eqref{pde_apriori}.
  Since $\bih{\cdot}{\cdot}$ is consistent, we have
  \begin{equation}\label{eq:duality_one}
    \begin{split}
      \Norm{R w}_{\leb{2}(\W)}^2
      &=
     -\int_{\W} \Delta g Rw \dx
      =
      -\int_{\W} \Delta g w \dx -\int_{\W} \Delta g (Rw-w) \dx\\
     & =
      \int_{\W} Rw w \dx +\int_{\W} \nabla g \cdot\nabla_h (Rw-w) \dx-\int_{\skel} \avg{\nabla g} \cdot\jump{Rw} \dS.
    \end{split}
  \end{equation}
  Let $\Pi: \sobh{1}(\W)\to \fes\cap \sobh{1}_0(\W)$ be a suitable conforming projection with optimal approximation properties, \highlight{for example the Cl\'ement interpolant}. Then, from the elliptic projection definition, we have
  \begin{equation}\label{eq:ell_orth_one}
  \int_{\W}\nabla_h Rw\cdot \nabla \Pi g\dx -\int_{\skel}\avg{\nabla \Pi g}\cdot\jump{Rw}\dS = \cA_h(Rw,\Pi g)
  =\cA_h(w,\Pi g)=\int_{\W}\nabla w\cdot \nabla \Pi g\dx. 
    \end{equation}
    Combining \eqref{eq:duality_one} with \eqref{eq:ell_orth_one}, we arrive at
    \begin{equation}\label{eq:duality_two}
    \begin{split}
      \Norm{R w}_{\leb{2}(\W)}^2     & =
      \int_{\W} Rw w \dx +\int_{\W} \nabla (g-\Pi g) \cdot\nabla_h (Rw-w) \dx-\int_{\skel} \avg{\nabla (g-\Pi g)} \cdot\jump{Rw} \dS\\
      &\le  \Norm{Rw}_{\leb{2}(\W)} \Norm{ w}_{\leb{2}(\W)} + \Norm{h^{-1}\nabla (g-\Pi g)}_{\leb{2}(\W)} \Norm{h \nabla_h (Rw-w)}_{\leb{2}(\W)}\\
      &\quad + \Norm{h^{-1/2}\avg{\nabla (g-\Pi g)}}_{\leb{2}(\skel)} \Norm{h^{1/2} \jump{Rw}}_{\leb{2}(\skel)}
     .
    \end{split}
  \end{equation}
From the optimal approximation properties of the projection/interpolant $\Pi$, we have
  \[
  \begin{split}
  &\Norm{h^{-1}\nabla_h \qp{ g - \Pi g}}_{\leb{2}(\W)}
  + \Norm{h^{-1/2}\avg{\nabla \qp{g - \Pi g}}}_{\leb{2}(\skel)}
  \le  C_{ap} \norm{g}_{\sobh{2}(\W)}\le C_{ap}C_{reg} \Norm{R w}_{\leb{2}(\W)}.
  \end{split}
  \]
  Setting $\widetilde{c}:= C_{ap}C_{reg}$ and combining the above, therefore, we arrive at
  \[
   \Norm{R w}_{\leb{2}(\W)} \le  \Norm{w}_{\leb{2}(\W)} +  \widetilde{c} \Big( \Norm{h \nabla_h \qp{Rw - w}}_{\leb{2}(\W)}
   + \Norm{h^{1/2}\jump{ R w}}_{\leb{2}(\skel)}\Big).
  \]
  It remains to show the bound
  \[
  \Norm{h \nabla_h Rw}_{\leb{2}(\W)}^2
   + \Norm{h^{1/2}\jump{ R w}}_{\leb{2}(\skel)}^2 
   \le C\big(  \Norm{h \nabla w}_{\leb{2}(\W)}^2
   + \Norm{h^{3/2}\avg{\nabla  w}}_{\leb{2}(\skel)}^2\big)+\frac{1}{4\widetilde c}\Norm{Rw}_{\leb{2}(\W)}^2,
  \]
  to conclude the proof. To that end, \eqref{eq:ip-proj} with $v_h=h^2R w$ implies
  \[
  \begin{split}
    \Norm{h \nabla_h Rw}_{\leb{2}(\W)}^2
    &=
    \int_{\W} \nabla_h R w \cdot \nabla_h \qp{h^2Rw}
      \dx
      \\
      &= 
      \int_{\W} \nabla  w \cdot \nabla_h \qp{h^2Rw}
      \dx
      \\
      &\qquad +
      \int_{\skel}\Big(
      \avg{\nabla (R w-w)}\cdot \jump{h^2R w} +\avg{\nabla (h^2 R w)}\cdot \jump{R w} 
      -
      \sigma \jump{R w}\cdot   \jump{h^2R w}\Big)\dS.
      \end{split}
  \]
Using now the elementary identities
$
 \jump{h^2R w}  =  \avg{h^2}\jump{R w}+ \jump{h^2}\avg{Rw}
$ 
 and $\jump{h^2} = 2\jump{h}\avg{h}$, which are valid on each internal face $e\in\skelint$, we arrive at
  \begin{equation}\label{mess}
  \begin{split}
    \Norm{h \nabla_h Rw}_{\leb{2}(\W)}^2  +  \int_{\skel} \sigma\avg{h^2} |\jump{R w}|^2\dS
    &= 
    \int_{\W} \qp{h\nabla  w} \cdot \qp{h\nabla_h Rw}
    \dx
    +
    \int_{\skel}
    \avg{\nabla (R w-w)}\cdot \jump{R w}\avg{h^2}\dS
    \\
    &\qquad
    +
    2\int_{\skelint} h \avg{\nabla (R w-w)}\cdot  \jump{h}\avg{Rw}\dS\\
    &
    \qquad
    + \int_{\skel}\avg{\nabla (h^2 R w)}\cdot \jump{R w} \dS
    -
    2\int_{\skelint}\sigma h\jump{R w}\cdot  \jump{h}\avg{Rw}\dS,       
  \end{split}
  \end{equation}
  recalling that $h:=\avg{h}$ on $\skelint$. Using \eqref{mesh_ass},  we proceed to bound each skeletal term on the right hand side of (\ref{mess}). To that end let  $C_{inv}>0$ denote the constant of the trace-inverse estimate, that is, $C_{inv}$ satisfies
  \[
  \Norm{v}_{\leb{2}(e)}^2\le C_{inv}k^2 h_K^{-1}\Norm{v}_{\leb{2}(K)}^2 \highlight{\quad \text{ for } v\in\fes},
  \]
  for $e\subset\partial K$, and $v\in \poly k (K)$, and recall $C_{qu}>0$ is the local quasi-uniformity constant from (\ref{eq:qu-const}). Then, in view of the definition of the penalty parameter (\ref{eq:penalty-param}), Cauchy-Schwarz and Young's inequalities, we see \begin{equation}
    \label{mess1}
\begin{split}
\int_{\skel}
      \avg{\nabla (R w-w)}\cdot \jump{R w} \avg{h^2}\dS
      \le &\
      \Norm{\sigma^{-1/2}\avg{h^2}^{1/2} \avg{\nabla (R w-w)}}_{\leb{2}(\skel)}^2 + \frac{1}{4}\Norm{\sqrt{\sigma\avg{h^2}}\jump{R w}}_{\leb{2}(\skel)}^2\\
     \le &\
     \frac{2C_{inv}C_{qu}}{k^2 C_{\sigma}}\Norm{h\nabla_h Rw}_{\leb{2}(\W)}^2 
     + \frac{2C_{qu}}{k^2 C_{\sigma}}\Norm{h^{3/2}\avg{\nabla w}}_{\leb{2}(\skel)}^2
     \\
     &\qquad
     + \frac{1}{4}\Norm{\sqrt{\sigma\avg{h^2}}\jump{R w}}_{\leb{2}(\skel)}^2,
\end{split}
  \end{equation}
  for the first skeletal term. Splitting the second up we have
  \begin{equation}
    \label{mess2}
    \begin{split}
      2  \int_{\skelint}
      \avg{\nabla R w}\cdot \jump{h}h\avg{Rw}\dS
\le &\ 
2  \alpha\int_{\skelint}
   h^2  | \avg{\nabla  \highlight{R} w}| |\avg{Rw}|\dS \\
\le &\  
\alpha  \Norm{h^{3/2}\avg{\nabla R w}}_{\leb{2}(\skelint)}^2
+ \alpha \Norm{h^{1/2}\avg{R w}}_{\leb{2}(\skelint)}^2\\
\le &\  
\alpha C_{inv}C_{qu} \big(\Norm{h\nabla_h R w}_{\leb{2}(\W)}^2
+ \Norm{R w}_{\leb{2}(\W)}^2\big).
    \end{split}
  \end{equation}
  Analogously we have
    \begin{equation}
    \label{mess3}
    \begin{split}
      2  \int_{\skelint}
      \avg{\nabla  w}\cdot \jump{h}h\avg{Rw}\dS
\le&\  
\alpha  \Norm{h^{3/2}\avg{\nabla  w}}_{\leb{2}(\skelint)}^2
+ \alpha C_{inv}C_{qu} \Norm{R w}_{\leb{2}(\W)}^2.
    \end{split}
    \end{equation}
    A similar argument to (\ref{mess1}) shows
    \begin{equation}
      \label{mess4}
  \begin{split} \int_{\skel}
      \avg{\nabla (h^2 R w)}\cdot \jump{Rw}\dS
\le &\
 \Norm{(\sigma\avg{h^2})^{-1/2}\avg{h^2\nabla R w}}_{\leb{2}(\skel)}^2
+ \frac{1}{4}\Norm{\sqrt{\sigma\avg{h^2}}\jump{R w}}_{\leb{2}(\skel)}^2\\
\le &\
\frac{C_{inv}C_{qu}}{k^2 C_{\sigma}} \Norm{h\nabla_h R w}_{\leb{2}(\W)}^2
+ \frac{1}{4}\Norm{\sqrt{\sigma\avg{h^2}}\jump{R w}}_{\leb{2}(\skel)}^2,
\end{split}
    \end{equation}
  and
  \begin{equation}
    \label{mess5}
    \begin{split}
      2\int_{\skelint}\kern -.2cm
      \sigma\avg{h}\jump{R w}\cdot  \jump{h}\avg{Rw}\dS
      \le&\ 2\alpha
      \int_{\skelint}
      \sigma h^2|\jump{R w}||\avg{Rw}|\dS
      \\
      \le &\ 
      4\alpha\Norm{\sigma^{1/2}h^2\avg{h^2}^{-1/2}\avg{R w}}_{\leb{2}(\skelint)}^2
      + \frac{\alpha}{4}\Norm{\sqrt{\sigma\avg{h^2}}\jump{R w}}_{\leb{2}(\skelint)}^2
      \\
      \le &\ 
      {4\alpha C_{inv}C_{qu} k^2 C_{\sigma}} \Norm{R w}_{\leb{2}(\W)}^2
      + \frac{\alpha}{4}\Norm{\sqrt{\sigma\avg{h^2}}\jump{R w}}_{\leb{2}(\skelint)}^2,
    \end{split}
  \end{equation}
  Substituting the above estimates (\ref{mess1})--(\ref{mess5}) into \eqref{mess}, we deduce 
  \[
  \begin{split}
    \Big(
    \frac{1}{2}
    -
    \alpha C_{inv}C_{qu}
    -
    \frac{3C_{inv}C_{qu}}{k^2 C_{\sigma}}
    \Big)
    \Norm{h \nabla_h Rw}_{\leb{2}(\W)}^2
    +
    &\frac{2-\alpha}{4}\int_{\skel} \sigma\avg{h^2} |\jump{R w}|^2\dS\\
    &\le  
    \frac{1}{2} \Norm{ h\nabla  w}_{\leb{2}(\W)}^2 
    +
    \qp{
      \alpha + \frac{2C_{qu}}{k^2 C_\sigma}
    }\Norm{h^{3/2}\avg{\nabla w}}_{\leb{2}(\skel)}^2
    \\
    &\qquad +
    2\alpha C_{inv}C_{qu}\big(1+2 k^2 C_{\sigma}\big)  \Norm{R  w}_{\leb{2}(\W)}^2 .
      \end{split}
  \]
  Therefore, assuming that the discontinuity-penalisation constant $C_{\sigma}$ is chosen so that 
  \begin{equation}
    \label{eq:penalty-param}
  C_{\sigma}>\max \{4, 24 k^{-2} C_{inv}C_{qu}, 8C_{inv}C_{qu}\},
  \end{equation}
   (with $C_{\sigma}>8C_{inv}C_{qu}$ being necessary for coercivity), upon selecting
  \begin{equation}\label{choice_of_alpha}
  \alpha <\min \Big\{1,  \frac{1}{8 C_{inv}C_{qu}}, \frac{\widetilde{c}}{32C_{inv}C_{qu}(2k^2C_{\sigma}+1)}\Big\},
  \end{equation}
  we have
   \[
  \begin{split}
    \frac{1}{4}\Norm{h \nabla_h Rw}_{\leb{2}(\W)}^2  +  \frac{2-\alpha}{4}\int_{\skel} \sigma\avg{h^2} |\jump{R w}|^2\dS
       &\le  
    \frac{1}{2} \Norm{ h\nabla  w}_{\leb{2}(\W)}^2 
      +
       C \Norm{h^{3/2}\avg{\nabla w}}_{\leb{2}(\skel)}^2
       \\
       &\qquad +
\frac{\widetilde{c}}{16}  \Norm{R  w}_{\leb{2}(\W)}^2 .
      \end{split}
  \]
  Notice that this is not the only choice of $C_\sigma$ and $\alpha$. There is a subtle dependency between the two values in that they are coupled such that choosing a larger $C_\sigma$ allows more flexibility on the selection of $\alpha$.
  
 The result already follows by combining the above bounds.
 \end{Proof}

\begin{Rem}[Polynomial degree dependence] \label{polydeg}
	\highlight{We remark on the dependence of the constants in \eqref{mesh_ass} on the polynomial degree $k$ in the proof of Theorem \ref{the:inf-sup} via its use of Lemma \ref{lem:bound-for-R}, as opposed to the more familiar (related) condition $ \|\nabla \widetilde h\|_{\leb{\infty}(\W)}\lesssim 1$ used in \cite{Makridakis:2016} for the proof of the respective result for the conforming finite element method. In \cite{Makridakis:2016}, the classical super-approximation argument given in \cite{NitscheSchatz:1974} is used, which is based on the repeated application of inverse estimates of the form $\|\nabla v\|_{\leb{2}(K)}\le Ck^2/h_K\|v\|_{K}$ for polynomials $v$ of degree $k$ in $K$ (cf., also \cite{Eriksson:1994,DemlowStevenson:2011}. Therefore, the respective bound on $ \|\nabla \widetilde h\|_{\leb{\infty}(\W)}\lesssim 1$ that is required to be small enough is proportional to $k^{2k}$. In the present proof, however, we \emph{avoided} the use of such super-approximation arguments in the proof of Lemma \ref{lem:bound-for-R}. As a result, the dependence of (small enough) constant required to satisfy the condition $\Norm{ \jump{h}/\avg{h}}_{\leb{\infty}(\Gamma)}$ is \emph{only} dependent on $k^2$, as seen by \eqref{choice_of_alpha}.  In fact, the mesh condition \eqref{mesh_ass} may be replaced by the condition
		\begin{equation}\label{mess_ass_2}
		\Norm{k^2 \jump{h}/\avg{h}}_{\leb{\infty}(\skelint)}\lesssim 1,
		\end{equation}
		with the right-hand side of the above inequality being \emph{independent} on the polynomial degree of the finite element space. We refrain from using the latter version of the mesh condition, however, in the interest of simplicity of the presentation.}
	
\highlight{The use of the condition \eqref{mess_ass_2} seems more natural to us in the dG setting. Although a very detailed comparison of the mesh conditions \eqref{connection} and \eqref{mess_ass_2} is beyond the scope of the present paper, it is clear that these conditions are qualitatively comparable. In addition,  the above argument indicates that in the context of high-order elements the approach taken herein might have an advantage over the use of a super-approximation argument.  Nonetheless, the numerical results in Figure 1 suggest that \eqref{mesh_ass} (or, equivalently, \eqref{mess_ass_2}) are reasonable in the context of non-globally quasiuniform (e.g., adaptive) meshes.}
        
\end{Rem}
 
 \begin{Rem}[Nonsymmetric interior penalty methods]
\label{rem:ns-int-pen}
   For $\theta\in [-1,1)$,   \eqref{eq:duality_two} becomes
\begin{equation}
    \begin{split}
      \Norm{R w}_{\leb{2}(\W)}^2     & =
      \int_{\W} Rw w \dx +\int_{\W} \nabla (g-\Pi g) \cdot\nabla_h (Rw-w) \dx\\
      &\quad-\int_{\skel} \avg{\nabla (g-\Pi g)} \cdot\jump{Rw} \dS-\int_{\skel}(1-\theta) \avg{\nabla \Pi g} \cdot\jump{Rw} \dS,
    \end{split}
  \end{equation}
 which, in turn, implies
   \[
   \Norm{R w}_{\leb{2}(\W)} \le  \Norm{w}_{\leb{2}(\W)} +  \widetilde{c} \Big( \Norm{h \nabla_h \qp{Rw - w}}_{\leb{2}(\W)}
   + \Norm{h^{1/2}\jump{ R w}}_{\leb{2}(\skel)}\Big) -(1-\theta)\frac{\int_{\skel} \avg{\nabla \Pi g} \cdot\jump{Rw} \dS}{ \Norm{R w}_{\leb{2}(\W)} }. 
  \]
Estimating the last term on the right-hand side of the above bound gives
  \[
  \begin{split}
  \int_{\skel} \avg{\nabla \Pi g} \cdot\jump{Rw} \dS
  \le&\  \Norm{(\sigma\avg{h^2})^{-1/2}\avg{\nabla (g-\Pi g}}_{\leb{2}(\skel)}\Norm{\sqrt{\sigma\avg{h^2}}\jump{Rw}}_{\leb{2}(\skel)}\\
  &+  \Norm{(\sigma\avg{h^2})^{-1/2}\avg{\nabla g}}_{\leb{2}(\skel)}\Norm{\sqrt{\sigma\avg{h^2}}\jump{Rw}}_{\leb{2}(\skel)}
  \\
  \le &\  \frac{2C_{reg}}{k \sqrt{\min\{h\}C_{\sigma}}}\Norm{Rw}_{\leb{2}(\W)}\Norm{\sqrt{\sigma\avg{h^2}}\jump{Rw}}_{\leb{2}(\skel)},
  \end{split}
  \]
  which yields \eqref{key} only with values of $\theta>1-\min\{h\}$ when $C_{\sigma}$ is chosen independent of $h$. When $C_{\sigma}$ is chosen to depend on negative powers of $h$, i.e., in the case of super-penalisation, we can retrieve \eqref{key} for non-symmetric versions of the interior penalty dG method.
 \end{Rem}

\subsection{Proof of Theorem \ref{the:inf-sup}} We give the proof of Theorem \ref{the:inf-sup} for $\theta=1$.
Our goal is for fixed $v_h\in\fes$ to construct a $w_h \in \fes$ such that
\begin{equation}
  \bih{w_h}{v_h} \geq \Norm{v_h}_{\leb{2}(\W)}^2
\end{equation}
and then to show one can find a constant $C>0$ such that
\begin{equation}
  \eenorm{w_h} \leq C \Norm{v_h}_{\leb{2}(\W)}.
\end{equation}
It is the case that each of the four components of the $\eenorm{\cdot}$-norm must be controlled; we shall bound each \highlight{of these} separately.

\subsubsection{Step 1:}
\label{step1}

For fixed $v_h$, let $\Phi\in\fes$ be the solution of the dual problem
\begin{equation}
  \bih{\Psi}{\Phi} = \ltwop{v_h}{\Psi} \Foreach \Psi \in \fes.
\end{equation}
To control $\enorm{\cdot}$, coercivity \eqref{eq:coer} yields
\begin{equation}
  \begin{split}
    \enorm{\Phi}^2
    \leq
    \frac 1 {c_0} \bih{\Phi}{\Phi}
    =
    \frac 1 {c_0} \ltwop{v_h}{\Phi}
    \leq
    \frac 1 {c_0} \Norm{v_h}_{\leb{2}(\W)}\Norm{\Phi}_{\leb{2}(\W)}
    \leq
    \frac {C_P}{c_0} \Norm{v_h}_{\leb{2}(\W)}\enorm{\Phi},
  \end{split}
\end{equation}
through a discrete Poincar\'e inequality and hence
\begin{equation}
  \enorm{\Phi} \leq \frac{C_P}{c_0} \Norm{v_h}_{\leb{2}(\W)}.
\end{equation}

\subsubsection{Step 2:}
\label{step2}

Let $w_1\highlight{\vert_K} = -b_K^2 \Delta \Phi$, for with $b_K$ denoting the standard polynomial bubble function vanishing on $\partial K$. We, then, have
\begin{equation}
  \begin{split}
    \bih{w_1}{\Phi}
= 
    \sum_{K\in\T{}}
    \int_K
    -\Delta \Phi w_1 \dx
    =\sum_{K\in\T{}}
  \int_K
  \norm{\Delta \Phi}^2 b_K^2
    \dx,
  \end{split}
\end{equation}
since $w_1\vert_e = 0$ and $\nabla w_1\vert_e = \vec 0$ for all $e\in\E$. 
Due to the equivalence of norms on finite dimensional linear spaces
\begin{equation}
  \label{eq:step2-1}
  \begin{split}
    \frac 1 {C_1}
    \sum_{K\in\T{}}
    \Norm{\Delta \Phi}^2_{\leb{2}(K)}
    &\le
    \int_K
    \norm{\Delta \Phi}^2 b_K^2
    \dx =   \bih{w_1}{\Phi}
    =
    \bih{\Phi}{w_1}   
    \\
    =&\ \bih{\Phi}{R w_1}
    =\bih{R w_1}{\Phi}
   =
    \ltwop{v_h}{R w_1}
   \leq
    \Norm{v_h}_{\leb{2}(\W)}\Norm{R w_1}_{\leb{2}(\W)},
  \end{split}
\end{equation}
using the symmetry of the bilinear form $\cA_h$. Recalling that $w_1$ is discrete, making use of Lemma \ref{lem:bound-for-R} and of local inverse inequalities, we have 
\begin{equation}
  \label{eq:step2-2}
  \begin{split}
    \Norm{R w_1}_{\leb{2}(\W)}
    &\leq 
    C \qp{\Norm{h \nabla w_1}_{\leb{2}(\W)} + \Norm{w_1}_{\leb{2}(\W)}
      + \Norm{h^{3/2}\avg{\nabla w_1}}_{\leb{2}(\skel)}
    }
    \\
    &\leq
    C \Norm{w_1}_{\leb{2}(\W)}
    \leq
    C \Big(\sum_{K\in\T{}} \Norm{\Delta \Phi}^2_{\leb{2}(K)}\Big)^{1/2}.
  \end{split}
\end{equation}
Combining (\ref{eq:step2-1}) and (\ref{eq:step2-2}), we see
\begin{equation}\label{eq:step_one}
  \Big(\sum_{K\in\T{}} \Norm{\Delta \Phi}^2_{\leb{2}(K)}\Big)^{1/2}
  \leq
 C \Norm{v_h}_{\leb{2}(\W)}.
\end{equation}

\subsubsection{Step 3:}
\label{step3}

Let $e$ be an internal edge of two neighbouring elements $K$ and $K'$ and let $b_e$ be a polynomial bubble function vanishing on the boundary of $K\cup K'\cup e$, so as to have by construction that $\nabla b_e\cdot \geovec n_e |_e=0$; the simplest such bubble function is of degree four when $d=2$, is a bubble function on the largest rhombus $\widetilde{K}_e$ contained fully in $K\cup K'\cup e$ and having $e$ as one of its diagonals $e$ (see, e.g., \cite{GeorgoulisHoustonVirtanen:2011} for details of such a construction). A completely analogous construction when $d=3$ yields the same properties. The mesh regularity assumed implies that $\diam (\widetilde{K}_e)$ is uniformly bounded above and below by the mesh-function $h$.

Let also $v_e:K\cup K'\cup e\to\mathbb{R}$ given by $v_e :=h^{-1}\jump{\nabla\Phi}$ on the face $e$ and extended as a constant on the direction of the normal to $e$. Setting 
\[
w_2:= \sum_{e\in \E} b_e^2v_e,
\]
we have $w_2\in \hoz$ and that $\nabla w_2\cdot \geovec n_e |_e=0$ for all $e\in\E$. Therefore,
\begin{equation}
  \begin{split}
    \bih{w_2}{\Phi}
    &=
    \sum_{K\in\T{}}
    \int_K
    -\Delta \Phi w_2
    \dx
    +
    \int_{\skelint}
    \jump{\nabla \Phi} \cdot \avg{w_2}
    \dS
    \\
    & =
  \sum_{K\in\T{}}
  \int_K
  -\Delta \Phi w_2
  \dx
  +
  \int_{\skelint}h^{-1}
  \norm{\jump{\nabla \Phi}}^2 b_e^2
  \dS.
  \end{split}
\end{equation}
The equivalence of all norms of a finite dimensional linear space implies that there exists a constant $C_2>0$, independent of $\Phi$ and $h$, such that
\begin{equation}
  \begin{split}
    \frac{1}{C_2} \int_{\skelint} h^{-1}\norm{\jump{\nabla \Phi}}^2 \dS
    &\leq
    \int_{\skelint}h^{-1} \norm{\jump{\nabla \Phi}}^2 b_e^2 \dS
    \leq
    \bih{w_2}{\Phi}
    +
    \sum_{K\in\T{}} \int_K \Delta \Phi w_2 \dx
    \\
    &=\
    \bih{R w_2}{\Phi}
    +
    \sum_{K\in\T{}} \int_K \Delta \Phi w_2 \dx
    =
    \ltwop{v_h}{R w_2}
    +
    \sum_{K\in\T{}} \int_K \Delta \Phi w_2 \dx
    \\
    &\leq
    C \Norm{v_h}_{\leb{2}(\W)} \qp{\Norm{R w_2}_{\leb{2}(\W)} + \Norm{w_2}_{\leb{2}(\W)}
      +\Norm{h^{3/2}\avg{\nabla w_2}}_{\leb{2}(\skel)}}
    \\
    &\leq C \Norm{v_h}_{\leb{2}(\W)}  \Norm{w_2}_{\leb{2}(\W)} ,
  \end{split}
\end{equation}
making use of \eqref{eq:step_one}, \eqref{key} and of standard inverse estimates, respectively. 
To finish, we observe the bound 
\[
\Norm{w_2}_{\leb{2}(\W)} \le C\Norm{\sqrt{h} w_2}_{\leb{2}(\skelint)} \le 
C\Norm{h^{-1/2}\jump{\nabla \Phi}}_{\leb{2}(\skelint)}, 
\] 
which, in turn, implies
\begin{equation}
    \Norm{h^{-1/2}\jump{\nabla \Phi}}_{\leb{2}(\skelint)}\le  C\Norm{v_h}_{\leb{2}(\W)} .
\end{equation}

\subsubsection{Step 4:}
\label{step4} As before, let $e$ be an internal edge of two neighbouring elements $K$ and $K'$ and let $b_e$ as in Step 3. Let also $p^{\ell}_e:K\cup K'\cup e\to \mathbb{R}$ be the plane passing through $e$ with slope equal to $h^{-3}$. Then, upon defining the function $z_e|_e:=(\Phi|_{\partial K\cap e}-\Phi|_{\partial K'\cap e})$ extended as a constant in the direction normal to $e$,  we set $w_3:\hoz\to \mathbb{R}$ given by
\[
w_3:= \sum_{e\in \E}z_eb_e^2 p_e^{\ell},
\]
where $z_e|_e:=\Phi|_{\partial K\cap e}$ is on the boundary faces $e\subset \partial\W$.
We note that the sign of the jump of $\Phi$ in the definition of $z_e$ is of no significance in what follows, so no effort is made in determining it exactly. With these definitions, we have $w_3=0$ on $\skelint$ and  
\[
\nabla w_3\cdot\geovec n_e |_e= h^{-3}\jump{\Phi}|_eb_e^2,
\]
on each $e\in\E$.
%
Therefore, we have 
\begin{equation}
  \begin{split}
    \bih{w_3}{\Phi}
    &=
    \sum_{K\in\T{}}
    \int_K
    -\Delta \Phi w_3
    \dx
    +\int_{\skel}
    \frac{b_e^2}{h^3} \norm{\jump{\Phi}}^2
    \dS.
  \end{split}
\end{equation}
As before, there exists a constant $C_2>0$, independent of $\Phi$ and of $h$, such that
\begin{equation}
  \begin{split}
    \frac 1 {C_2}\int_{\skel}
    h^{-3} \norm{\jump{\Phi}}^2
    \dS
    &\le
    \int_{\skel}
    \frac{b_e^2}{h^3} \norm{\jump{\Phi}}^2
    \dS
    =
    \bih{w_3}{\Phi}
    +
    \sum_{K\in\T{}}
    \int_K
    \Delta \Phi w_3
    \dx\\
   & \leq
    \Norm{v_h}_{\leb{2}(\W)} \Norm{R w_3}_{\leb{2}(\W)}
    +
    \Big(\sum_{K\in\T{}} \int_K \norm{\Delta \Phi}^2 \dx\Big)^{1/2} \Norm{w_3}_{\leb{2}(\W)}\\
    &    \le  C \Norm{v_h}_{\leb{2}(\W)} \qp{\Norm{R w_3}_{\leb{2}(\W)} + \Norm{w_3}_{\leb{2}(\W)}
      +\Norm{h^{3/2}\avg{\nabla w_3}}_{\leb{2}(\skel)}}
    \\
    &\le C \Norm{v_h}_{\leb{2}(\W)}  \Norm{w_3}_{\leb{2}(\W)},
  \end{split}
\end{equation}
from \eqref{eq:step_one}, \eqref{key} and standard inverse estimates.  

Also, we have
\[
 \Norm{w_3}_{\leb{2}(\W)}^2
 \le 
 C\sum_{K\in\T{},e\subset\partial K}\Norm{p_e^{\ell}}_{\leb{\infty}(K)}^2 \Norm{z_e}_{\leb{2}(K)}^2
 \leq
C\sum_{e\subset\partial \E}h_K^{-4} \Norm{\sqrt{h}\jump{\Phi}}_{\leb{2}(e)}^2
\leq
C \Norm{h^{-3/2}\jump{\Phi}}_{\leb{2}(\skel)}^2,
\]
which, finally, implies
\begin{equation}
\Norm{h^{-3/2}\jump{\Phi}}_{\leb{2}(\skel)} \le C \Norm{v_h}_{\leb{2}(\W)},
\end{equation} 
which, \highlight{taking $w_h = R\qp{w_1+w_2+w_3}$}, already proves the result.
%

\section{Relaxation of regularity requirements}

In the above discussion, we assumed for clarity of presentation that
for the exact solution we have $u\in\sobh{2}(\W)$; the analysis presented also holds if
$u\in\sobh{s}(\W)$ for $s > 3/2$. In this section we shall deduce a useful a priori bound for the interior penalty method with $\theta =1$ for the case
$u\in\sobh{1}(\W)$ also, by showing that 
\begin{equation}
  \Norm{u - u_h}_{\leb{2}(\W)}
  \leq 
  C\qp{
    \inf_{w_h\in\fes} \Norm{u - w_h}_{\leb{2}(\W)}
    +
      \Norm{h^{2}(f - P_{k} f)}_{\leb{2}(\W)}},
\end{equation}
where $P_k$ is the $\leb{2}$-orthogonal projection operator into element-wise polynomials of degree $k$. To do so, we shall use ideas from \cite{Gudi:2010}, extended to the present setting through the following Lemmata. The main result of the section is stated in Theorem \ref{the:optimal-weak}. 

\begin{Lem}[]
  \label{lem:1}
  For $w \in \sobh{1}(\W)$, $v\in\sobh{2}(\W)$ and $w_h\in\fes$ it holds that
  \begin{equation}
    \norm{ 
      \bi{w}{v} 
      - 
      \bih{w_h}{v}
    }
    \leq
    C \Norm{w - w_h}_{\leb{2}(\W)} \eenorm{v}.
  \end{equation}
\end{Lem}

\begin{Proof}
  Since $v\in\sobh{2}(\W)$ we have, through the consistency of the scheme, that
  \begin{equation}
    \norm{ 
      \bi{w}{v} 
      - 
      \bih{w_h}{v}
    }
    =
    \norm{
      - \int_\W \qp{w - w_h} \Delta v
    }
    \leq
    \Norm{w - w_h}_{\leb{2}(\W)} \eenorm{v}.
  \end{equation}
\end{Proof}

\begin{Lem}[Reconstruction operator]
  \label{lem:2}
  Let $\HCT$ denote the Hsieh-Clough-Tocher macro-element space, then there exists an operator $E : \fes \to \HCT \subset \sobh{2}(\W)$ such that for $\alpha = 0,1,2$
  \begin{equation}
    \sum_{K\in\T{}} \Norm{E (w_h) - w_h}_{\sobh{\alpha}(K)}^2
    \leq
    C \qp{\Norm{h^{1/2-\alpha}\jump{w_h}}_{\leb{2}(\Gamma)}^2
      +
      \Norm{h^{3/2-\alpha}\jump{\nabla w_h}}_{\leb{2}(\Gamma)}^2} \Foreach w_h\in\fes.
  \end{equation}
\end{Lem}
\begin{Proof}
  The proof of this is given in \cite[Lemma 3.1]{GeorgoulisHoustonVirtanen:2011}.
\end{Proof}

\begin{Rem}[]
  The use of Lemma \ref{lem:1} will be crucial subsequently and, for this reason, the use of an $\sobh{2}$ reconstruction operator is necessary. This is in contrast to the argument in \cite{Gudi:2010} where an $\sobh{1}$ conforming reconstruction was used. 
\end{Rem}

\begin{Lem}[A posteriori lower bound]
  \label{lem:3}
  Let $u\in\sobh{s}(\W)$ be the weak solution of (\ref{eq:weakform}) and $w_h \in\fes$ be an arbitrary finite element function. Then,
    \begin{equation}
    \sup_{v_h \in \fes}\frac{\ltwop{f}{v_h - E(v_h)} - \bih{w_h}{v_h - E(v_h)}}{\eenorm{v_h}}
    \leq
    C \qp{
      \sum_{K\in\T{}}
      \Norm{u - w_h}_{\leb{2}(K)}^2
      +
        \Norm{h^{2}(f - P_{k} f)}_{\leb{2}(K)}^2
      }^{1/2}.
  \end{equation}
\end{Lem}
\begin{Proof}
  We begin by noting that
  \begin{equation}
    \begin{split}
      \ltwop{f}{v_h - E(v_h)} - \bih{w_h}{v_h - E(v_h)}
      &=
      \sum_{K\in\T{}}\int_K \qp{f + \Delta w_h}\qp{v_h - E(v_h)}\dx
      \\
      &\qquad
      - \int_{\skelint} \jump{\nabla w_h} \avg{v_h - E(v_h)}\dS
      \\
      &\qquad
      +
      \int_{\skel}\big( \jump{w_h}\cdot \avg{\nabla v_h - \nabla E(v_h)}
      - \sigma \jump{w_h} \cdot \jump{v_h}\big)\dS
      \\
      &=: \sum_{K\in\T{}} \cI_{1,K} + \sum_{e\in\E} \cI_{2,e} + \cI_{3,e} + \cI_{4,e}.
    \end{split}
  \end{equation}
  We proceed to control each term separately. Firstly,
  \begin{equation}
    \label{eq:pf0}
    \begin{split}
      \cI_{1,K}
      &\leq
      \Norm{h^2\qp{f + \Delta w_h}}_{\leb{2}(K)} \Norm{h^{-2} \qp{v_h - E(v_h)}}_{\leb{2}(K)}
      \\
      &\leq
      \qp{\Norm{h^2\qp{P_{k} f + \Delta w_h}}_{\leb{2}(K)} + \Norm{h^2\qp{f - P_{k} f}}_{\leb{2}(K)}} \Norm{h^{-2} \qp{v_h - E(v_h)}}_{\leb{2}(K)}.
    \end{split}
  \end{equation}
  Now, making use of the interior bubble function $b_K$,
  we have
  \begin{equation}
    \label{eq:pf1}
    \begin{split}
      \Norm{h^2\qp{P_{k} f + \Delta w_h}}_{\leb{2}(K)}^2
      &\leq
      C_1 \int_K
      h^4 \qp{P_{k} f + \Delta w_h} b^2_K \qp{P_{k} f + \Delta w_h}\dx
      \\
      &=
      C_1 \int_K h^4\qp{\qp{P_{k} f  - f} + \qp{f + \Delta w_h}} b^2_K \qp{P_{k} f + \Delta w_h}\dx
      \\
      &\leq
      C_1 \Norm{h^2\qp{P_{k} f  - f}}_{\leb{2}(K)}
      \Norm{h^2 b_K^2 \qp{P_{k} f + \Delta w_h}}_{\leb{2}(K)} \\
      &\quad
      + C_1\int_K h^4\qp{\qp{f + \Delta w_h}} {b^2_K \qp{P_{k} f + \Delta w_h}}\dx.
    \end{split}
  \end{equation}
  Since $b^2_K = 0$ and $\nabla b^2_K = \vec 0$ on $\partial K$, we have
  \begin{equation}
    \label{eq:pf2}
    \begin{split}
      \int_K h^4\qp{f + \Delta w_h} {b^2_K \qp{P_{k} f + \Delta w_h}}\dx
      &=
      \int_K h^4\qp{u - w_h} \cdot \Delta \qp{b_K^2 \qp{P_{k} f + \Delta w_h}}\dx
      \\
      &\leq
     C\Norm{u - w_h}_{\leb{2}(K)} \Norm{h^2\qp{P_{k} f + \Delta w_h}}_{\leb{2}(K)},
    \end{split}
  \end{equation}
  making use of inverse inequalities.
  Hence combining (\ref{eq:pf0}) and (\ref{eq:pf1}) with (\ref{eq:pf2}) we see
  \begin{equation}
    \label{eq:pf3}
    \begin{split}
      \cI_{1,K}
      &\leq
      C
      \qp{\Norm{u - w_h}_{\leb{2}(K)} + \Norm{h^2\qp{P_{k} f  - f}}_{\leb{2}(K)}}
      \Norm{h^{-2} \qp{v_h - E(v_h)}}_{\leb{2}(K)}.
    \end{split}
  \end{equation}

  Secondly,
  \begin{equation}
    \begin{split}
      \cI_{2,e}
      &\leq
      \Norm{h^{3/2} \jump{\nabla w_h}}_{\leb{2}(e)}
      \Norm{h^{-3/2} \avg{v_h - E(v_h)}}_{\leb{2}(e)}.
    \end{split}
  \end{equation}
  Now
  \begin{equation}
    \begin{split}
      \Norm{h^{3/2} \jump{\nabla w_h}}^2_{\leb{2}(e)}
      &\leq
      C \int_e h^3 \jump{\nabla w_h} b^2_e \jump{\nabla w_h} \dS
      \\
      &\leq
      C \int_e h^4 \jump{\nabla w_h - \nabla u} b^2_e v_e \dS,
    \end{split}
  \end{equation}
  with $v_e$ defined in Step 3 of the Proof of Theorem \ref{the:inf-sup}. Now since $v_e b^2_e$ vanishes over the $\partial\qp{K\cup K'}$ and $\nabla b_e^2 \cdot \vec n = 0$ we see
  \begin{equation}
    \begin{split}
      \Norm{h^{3/2} \jump{\nabla w_h}}^2_{\leb{2}(e)}
      &\leq
      C \int_{K \cup K'}
      h^4\qp{\qp{ u - w_h} \Delta_h \qp{b_e^2 v_e}
      -
      \qp{ f + \Delta_h w_h} b_e^2 v_e}
      \d x
      \\
      &\leq
      C\qp{\Norm{u - w_h}_{\leb{2}(K \cup K')}
        +
        \Norm{h^2 \qp{f + \Delta_h w_h}}_{\leb{2}(K \cup K')}
      }
      \Norm{h^2 v_e b_e^2}_{\leb{2}(K \cup K')},
    \end{split}
  \end{equation}
  though inverse inequalities. Now note that in view of the properties of $b_e$ there exists a constant such that
  \begin{equation}
    \Norm{h^2 v_e b_e^2}_{\leb{2}(K \cup K')}
    \leq
    C
    \Norm{h^{3/2} \jump{\nabla w_h}}_{\leb{2}(e)}
  \end{equation}
  to see
  \begin{equation}
    \label{eq:pf4}
    \cI_{2,e} \leq C
    \qp{
      \Norm{u - w_h}_{\leb{2}(K \cup K')}
        +
        \Norm{h^2 \qp{f + \Delta_h w_h}}_{\leb{2}(K \cup K')}
    }
    \Norm{h^{-3/2} \avg{v_h - E(v_h)}}_{\leb{2}(e)}.
  \end{equation}
  The third term
  \begin{equation}
    \label{eq:pf5}
    \begin{split}
      \cI_{3,e}
      &\leq
      \Norm{h^{1/2} \jump{w_h - u}}_{\leb{2}(e)}
      \Norm{h^{-1/2} \avg{\nabla v_h - \nabla E(v_h)}}_{\leb{2}(e)}
      \\
      &\leq
      C_{inv}C_{qu}\Norm{w_h - u}_{\leb{2}(K\cup K')}
      \Norm{h^{-1/2} \avg{\nabla v_h - \nabla E(v_h)}}_{\leb{2}(e)}.
    \end{split}
  \end{equation}
  Finally the fourth term,
  \begin{equation}
    \label{eq:pf6}
    \begin{split}
      \cI_{4,e}
      &\leq
      C_\sigma k^2
      \Norm{h^{1/2} \jump{w_h - u}}_{\leb{2}(e)}
      \Norm{h^{-3/2} \jump{v_h - E(v_h)}}_{\leb{2}(e)}
      \\
      &\leq
      C_\sigma C_{inv} C_{qu} k^2
      \Norm{w_h - u}_{\leb{2}(e)}
      \Norm{h^{-3/2} \jump{v_h - E(v_h)}}_{\leb{2}(e)}.
    \end{split}
  \end{equation}
  Collecting all the information thusfar from (\ref{eq:pf3}), (\ref{eq:pf4}), (\ref{eq:pf5}) and (\ref{eq:pf6}) we can conclude that 
  \begin{equation}
    \begin{split}
      \ltwop{f}{v_h - E(v_h)} - \bih{w_h}{v_h - E(v_h)}
      \leq
      C
      \sum_{K\in\T{}}
      \qb{
        \qp{
          \Norm{u - w_h}_{\leb{2}(K)} + \Norm{h^2\qp{P_{k} f - f}}_{\leb{2}(K)}
        }
        \eta_K(v_h - E(v_h))
      },
    \end{split}
  \end{equation}
  where
  \begin{equation}
    \eta_K(z)
    =
    \max\qp{
      \Norm{h^{-2} z}_{\leb{2}(K)},
      \max_{e\in\partial K}
      \Norm{h^{-3/2} z}_{\leb{2}(e)},
      \max_{e\in\partial K}
      \Norm{h^{-1/2} \nabla z}_{\leb{2}(e)}
    }.
  \end{equation}
  Using the approximability properties of $E$ given in Lemma \ref{lem:2} we see
  \begin{equation}
    \sum_{K\in\T{}}
    \eta_K(v_h - E(v_h))
    \leq
    C
    \eenorm{v_h}^2,
  \end{equation}
  and hence the result follows from a discrete Cauchy-Schwarz inequality.
\end{Proof}

\begin{The}[Optimal convergence for weak solutions]
  \label{the:optimal-weak}
  Let $u\in\sobh{s}(\W)$ be the weak solution of
  (\ref{eq:bilinear-form}) and $u_h\in\fes$ be the dG solution (\ref{eq:IP}) and that the conditions of Theorem \ref{the:inf-sup} hold. Then there exists a constant such that
  \begin{equation}
    \Norm{u - u_h}_{\leb{2}(\W)}
    \leq
    C\qp{
      \inf_{w_h\in\fes} \Norm{u - w_h}_{\leb{2}(\W)}
      +
      \qp{
        \sum_{K\in\T{}} \Norm{h^2\qp{f - P_{k} f}}_{\leb{2}(K)}^2
      }^{1/2}
    },
  \end{equation}
  where $P_k$ is the $\leb{2}$ orthogonal projector into piecewise polynomials of degree $k$. 
\end{The}

\begin{Proof}[of Theorem \ref{the:optimal-weak}]
  Note that from Theorem \ref{the:inf-sup} we have that for any $w_h \in \fes$
  \begin{equation}
    \znorm{u_h - w_h}
    \leq
    \sup_{v_h \in \fes}\frac{\bih{u_h - w_h}{v_h}}{\eenorm{v_h}}.
  \end{equation}
  By adding and substracting appropriate terms we see
  \begin{equation}
    \bih{u_h - w_h}{v_h}
    =
    \bi{u}{E(v_h)}
    -
    \bih{w_h}{E(v_h)}
    +
    \ltwop{f}{v_h - E(v_h)}
    -
    \bih{w_h}{v_h - E(v_h)}
  \end{equation}
  and by Lemma \ref{lem:1}
  \begin{equation}
    \begin{split}
      \bi{u}{E(v_h)} - \bih{w_h}{E(v_h)}
      &\leq
      C
      \Norm{u - w_h}_{\leb{2}(\W)}
      \eenorm{E(v_h)}
      \\
      &\leq
      C\Norm{u - w_h}_{\leb{2}(\W)}
      \eenorm{v_h},
    \end{split}
  \end{equation}
  by Lemma \ref{lem:2}. Hence
  \begin{equation}
    \znorm{u_h - w_h}
    \leq
    C\qp{
      \Norm{u - w_h}_{\leb{2}(\W)}
      +
      \sup_{v_h \in \fes}\frac{\ltwop{f}{v_h - E(v_h)} - \bih{w_h}{v_h - E(v_h)}}{\eenorm{v_h}}
    }
    .
  \end{equation}
  So clearly,
  \begin{equation}
    \begin{split}
      \Norm{u - u_h}_{\leb{2}(\W)}
      &\leq
      \Norm{u - w_h}_{\leb{2}(\W)}
      +
      \Norm{u_h - w_h}_{\leb{2}(\W)}
      \\
      &\leq
      C\qp{\Norm{u - w_h}_{\leb{2}(\W)}
      +
      \znorm{u_h - w_h}}
      \\
      &\leq
      C\qp{
        \Norm{u - w_h}_{\leb{2}(\W)}
        +
        \sup_{v_h \in \fes}\frac{\ltwop{f}{v_h - E(v_h)} - \bih{w_h}{v_h - E(v_h)}}{\eenorm{v_h}}
      }
      .
    \end{split}
  \end{equation}
  The result follows from Lemma \ref{lem:3}.
\end{Proof}

\end{document}